\documentclass[11pt]{article}
\usepackage{graphicx}
\usepackage{color}
\usepackage{amsmath}
\usepackage{amssymb}
\usepackage{amscd}
\usepackage{setspace}

\newcommand{\R}{\mathbb{R}}
\newcommand{\inr}[1]{\bigl< #1 \bigr>}

\newcommand{\E}{\mathbb{E}}

\newcommand{\eps}{\varepsilon}

\newcommand \cL{{\cal L}}

\newcommand \bA{{\mathbb A}}

\newcommand{\Pro}{\mathbb P}
\newcommand{\norm}[1]{\left\|#1\right\|}%

\newtheorem{Theorem}{Theorem}[section]
\newtheorem{Lemma}[Theorem]{Lemma}

\newtheorem{Definition}[Theorem]{Definition}

\newtheorem{Corollary}[Theorem]{Corollary}
\newtheorem{Remark}[Theorem]{Remark}

\newtheorem{Assumption}{Assumption}[section]

\numberwithin{equation}{section}

\def \proof {\noindent {\bf Proof.}\ \ }

\def \endproof
{{\mbox{}\nolinebreak\hfill\rule{2mm}{2mm}\par\medbreak}}

\begin{document}
\title{{Minimax rate of convergence and the performance of ERM in phase recovery}}
\author{Guillaume Lecu\'e${}^{1,3}$  \and Shahar Mendelson${}^{2,4,5}$}

\footnotetext[1]{CNRS, CMAP, Ecole Polytechnique, 91120 Palaiseau, France.}
\footnotetext[2]{Department of Mathematics, Technion, I.I.T, Haifa
32000, Israel.}
 \footnotetext[3] {Email:
guillaume.lecue@cmap.polytechnique.fr }
\footnotetext[4] {Email:
shahar@tx.technion.ac.il}
\footnotetext[5]{Supported by the Mathematical Sciences Institute -- The Australian National University.}



\maketitle

\begin{abstract}
We study the performance of Empirical Risk Minimization in noisy phase retrieval problems, indexed by subsets of $\R^n$ and relative to subgaussian sampling; that is, when the given data is $y_i=\inr{a_i,x_0}^2+w_i$ for a subgaussian random vector $a$, independent noise $w$ and a fixed but unknown $x_0$ that belongs to a given subset of $\R^n$.

We show that ERM produces $\hat{x}$ whose Euclidean distance to either $x_0$ or $-x_0$ depends on the gaussian mean-width of the indexing set and on the signal-to-noise ratio of the problem. The bound coincides with the one for linear regression when $\|x_0\|_2$ is of the order of a constant. In addition, we obtain a minimax lower bound for the problem and identify sets for which ERM is a minimax procedure. As examples, we study the class of $d$-sparse vectors in $\R^n$ and the unit ball in $\ell_1^n$.
\end{abstract}

\section{Introduction}
Phase retrieval has attracted much attention recently, as it has natural applications in areas that include X-ray crystallography, transmission electron microscopy and coherent diffractive imaging (see, for example, the discussion in \cite{EM12} and references therein).

In phase retrieval, one attempts to identify a vector $x_0$ that belongs to an arbitrary set $T \subset \R^n$ using noisy, quadratic measurements of $x_0$. The given data is a sample of cardinality $N$, $(a_i,y_i)_{i=1}^N$, for vectors $a_i \in \R^n$ and
\begin{equation}
\label{eq:meas}
y_i=|\inr{a_i,x_0}|^2+w_i,
\end{equation}
for a noise vector $(w_i)_{i=1}^N$.

Our aim is to investigate phase retrieval from a theoretical point of view, relative to a well behaved, random sampling method. To formulate the problem explicitly, let $\mu$ be an isotropic, $L$-subgaussian measure on $\R^n$ and set $a$ to be a random vector distributed according to $\mu$. Thus, for every $x \in \R^n$, $\E \inr{x,a}^2=\|x\|_2^2$ (isotropicity) and for every $u \geq 1$, $Pr(|\inr{x,a}| \geq Lu\|\inr{x,a}\|_{L_2}) \leq 2\exp(-u^2/2)$ ($L$-subgaussian).


Given a set $T \subset \R^n$ and a fixed, but unknown $x_0 \in T$, $y_i$ are the random noisy measurements of $x_0$: for a sample size $N$, $(a_i)_{i=1}^N$ are independent copies of $a$ and $(w_i)_{i=1}^N$ are independent copies of a mean-zero variable $w$ that are also independent of $(a_i)_{i=1}^N$.

Clearly, due to the nature of the given measurements, $x_0$ and $-x_0$ are indistinguishable, and the best that one can hope for is a procedure that produces $\hat{x} \in T$ that is close to one of the two points.

The goal here is to find such a procedure and identify the way in which the distance between $\hat{x}$ and either $x_0$ or $-x_0$ depends on the structure of $T$, the measure $\mu$ and the noise.


The procedure studied here is empirical risk minimization (ERM), which produces $\hat{x}$ that minimizes the empirical risk in $T$:
\begin{equation*}
P_N \ell_x=\frac{1}{N}\sum_{i=1}^N \big( \inr{a_i,x}^2 - y_i \big)^2.
\end{equation*} 
The loss is the standard squared loss functional, which, in this case,satisfies 
\begin{equation*}
\ell_x(a,y) = (f_x(a)-y)^2= (\inr{x,a}^2 -\inr{x_0,a}^2-w)^2=(\inr{x-x_0,a}\inr{x+x_0,a}-w)^2.
\end{equation*}

Comparing the empirical and actual structures on $T$ is a vital component in the analysis of ERM. In phase recovery, the centered empirical process that is at the heart of this approach is defined for any $x\in T$ by,
\begin{equation*}
P_N (\ell_x-\ell_{x_0}) =  \frac{1}{N} \sum_{i=1}^N \inr{x-x_0,a_i}^2\inr{x+x_0,a_i}^2
-  \frac{2}{N} \sum_{i=1}^N w_i \inr{x-x_0,a_i}\inr{x+x_0,a_i}.
\end{equation*}
Both the first and second components are difficult to handle directly, even when the underlying measure is subgaussian, because of the powers involved (an effective power of $4$ in the first component and of $3$ in the second one).

Therefore, rather than trying to employ the concentration of empirical means around the actual ones, which might not be sufficiently strong in this case, one uses a combination of a small-ball estimate for the `high order' part of the process, and a more standard deviation argument for the low-order component (see Section \ref{sec:squared-loss} and the formulation of Theorem A and Theorem B).

We assume that linear forms satisfy a certain small-ball estimate, and in particular, do not assign too much weight to small neighbourhoods of $0$.
\begin{Assumption} \label{ass:small-ball}
There is a constant $\kappa_0>0$ satisfying that for every $s,t \in \R^n$,
\begin{equation*}
\E|\inr{a,s} \inr{a,t}| \geq \kappa_0 \|s\|_{2}\|t\|_{2}.
\end{equation*}
\end{Assumption}

Assumption \ref{ass:small-ball} is not very restrictive and holds for many natural choices of random vectors in $\R^n$, like the gaussian measure or any isotropic log-concave measure on $\R^n$ (see, for example, the discussion in \cite{EM12}).

It is not surprising that the error rate of ERM depends on the structure of $T$, and because of the subgaussian nature of the  random measurement vector $a$, the natural parameter that captures the complexity of $T$ is the gaussian mean-width associated with normalizations of $T$.

\begin{Definition} \label{def:ell}
Let $G=(g_1,...,g_n)$ be the standard gaussian vector in $\R^n$. For $T \subset \R^n$, set
\begin{equation*}
\ell(T)=\E \sup_{t \in T} \Big|\sum_{i=1}^n g_i t_i\Big|.
\end{equation*}
\end{Definition}
The normalized sets in question are
\begin{align*}
T_{-,R} & =\left\{\frac{t-s}{\|t-s\|_2} \ : \ t,s \in T, \ \ R < \|t-s\|_2\|t+s\|_2 \right\},
\\
T_{+,R} & =\left\{\frac{t+s}{\|t+s\|_2} \ : \ t,s \in T, \ \ R < \|t-s\|_2\|t+s\|_2 \right\},
\end{align*}
which have been used in \cite{EM12}, or their `local' versions,
\begin{align*}
T_{-,R}(x_0) & =\left\{\frac{t-x_0}{\|t-x_0\|_2} \ : \ t \in T, \ \ R < \|t-x_0\|_2\|t+x_0\|_2 \right\},
\\
T_{+,R}(x_0) & =\left\{\frac{t+x_0}{\|t+x_0\|_2} \ : \ t \in T, \ \ R < \|t-x_0\|_2\|t+x_0\|_2 \right\}.
\end{align*}
The sets in question play a central role in the exclusion argument that is used in the analysis of ERM. Setting ${\cal L}_x = \ell_x - \ell_{x_0}$, the excess loss function associated with $\ell$ and $x \in T$, it is evident that $P_N {\cal L}_{\hat{x}} \leq 0$ (because ${\cal L}_{x_0} = 0$ is a possible competitor). If one can find an event of large probability and $R>0$ for which $P_N {\cal L}_x >0$ if $\|x-x_0\|_2 \|x+x_0\|_2 \geq R$, then on that event, $\|\hat{x}-x_0\|_2 \|\hat{x}+x_0\|_2 \leq R$, which is the estimate one is looking for.

The normalization allows one to study `relative fluctuations' of $P_N {\cal L}_x$ -- in particular, the way these fluctuations scale with $\|x-x_0\|_2 \|x+x_0\|_2$. This is achieved by considering empirical means of products of functions $\inr{u,\cdot} \inr{v,\cdot}$, for $u \in T_{+,R}(x_0)$ and $v \in T_{-,R}(x_0)$.

The obvious problem with the `local' sets $T_{+,R}(x_0)$ and $T_{-,R}(x_0)$ is that $x_0$ is not known. As a first attempt of bypassing this problem, one may use the `global' sets $T_{+,R}$ and $T_{-,R}$ instead, as had been done in \cite{EM12}.

Unfortunately, this global approach is not completely satisfactory. Roughly put, there are two types of subsets of $\R^n$ one is interested in, and that appear in applications. The first type consists of sets for which the `local complexity' is essentially the same everywhere, and the sets $T_{+,R}, T_{-,R}$ are not very different from the seemingly smaller $T_{+,R}(x_0)$, $T_{-,R}(x_0)$, regardless of $x_0$.

A typical example of such a set is $d$-sparse vectors -- a set consisting of all the vectors in $\R^n$ that are supported on at most $d$-coordinates. For every $x_0 \in T$ and $R>0$, the sets $T_{+,R}(x_0),T_{-,R}(x_0)$, and $T_{+,R}, T_{-,R}$ are contained in the subset of the sphere consisting of $2d$-sparse vectors, which is a relatively small set.

For this kind of set, the `global' approach, using $T_{+,R}$ and $T_{-,R}$, suffices, and the choice of the target $x_0$ does not really influence the rate in which $\norm{\hat x-x_0}_2\norm{\hat x+x_0}_2$ decays to $0$ with $N$.

In contrast, sets of the second type one would like to study, have vastly changing local complexity, with the typical example being a convex, centrally symmetric set (i.e. if $x \in T$ then $-x \in T$).

Consider, for example, the case $T=B_1^n$, the unit ball in $\ell_1^n$. It is not surprising that for small $R$, the sets $T_{+,R}(0)$ and $T_{-,R}(0)$ are very different from $T_{-,R}(e_1)$ and $T_{+,R}(e_1)$: the ones associated with the centre $0$ are the entire sphere, while for $e_1=(1,0,....,0)$, $T_{+,R}(e_1)$ and $T_{-,R}(e_1)$ consist of vectors that are well approximated by sparse vectors (whose support depend on $R$), and thus are rather small subsets of the sphere .

The situation that one encounters in $B_1^n$ is  generic for convex centrally-symmetric sets. The sets become locally `richer' the closer the centre is to $0$, and at $0$, for small enough $R$, $T_{+,R}(0)$ and $T_{-,R}(0)$ are the entire sphere. Since the sets $T_{+,R}$ and $T_{-,R}$ are blind to the location of the centre, and are, in fact, the union over all possible centres of the local sets, they are simply too big to be used in the analysis of ERM in convex sets. A correct estimate on the performance of ERM for such sets requires a more delicate local analysis and additional information on $\|x_0\|_2$. Moreover, the rate of convergence of ERM {\it truly depends} on $\|x_0\|_2$ in the phase recovery problem via the signal-to-noise ratio $\norm{x_0}_2/\sigma$.

We begin by formulating our results using the `global' sets $T_{+,R}$ and $T_{-,R}$. Let $T_+ = T_{+,0}$ and $T_-=T_{-,0}$, set
\begin{equation*}
E_R=\max\{\ell(T_{+,R}),\ell(T_{-,R})\}, \ \ \ E=\max\{\ell(T_+),\ell(T_-)\}
\end{equation*}
and observe that as nonempty subsets of the sphere $\ell(T_{-,r}), \ell(T_{+,r}) \geq \E |g| = \sqrt{2/\pi}$.

The first result presented here is that the error rate of ERM for the phase retrieval problem in $T$ depends on the fixed points
\begin{equation*}
r_{2}(\gamma) = \inf \left\{r>0 : E_r \leq  \gamma\sqrt{N}r \right\}
\end{equation*}
and
\begin{equation*}
r_0(Q)=\inf \left\{r>0 : E_r \leq  Q\sqrt{N} \right\},
\end{equation*}
for constants $\gamma$ and $Q$ that will be specified later.

Recall that the $\psi_2$-norm of a random variable $w$ is defined by $\|w\|_{\psi_2}=\inf\{c>0:\E\exp(w^2/c^2)\leq2\}$ and set $\sigma=\|w\|_{\psi_2}$.

\vskip0.7cm
\noindent{\bf Theorem A.}
For every $L>1$, $\kappa_0>0$ and $\beta>1$, there exist constants $c_0,c_1$ and $c_2$ that depend only on $L$, $\kappa_0$ and $\beta$ for which the following holds. Let $a$ and $w$ be as above and assume that $w$ has a finite $\psi_2$-norm. If $\ell$ is the squared loss and $\hat{x}$ is produced by ERM, then with probability at least
$$
1-2\exp(-c_0 \min\{\ell^2(T_{+,r_2^*}), \ell^2(T_{-,r_2^*})\}) - 2N^{-\beta+1},
$$
$$
\|\hat{x}-x_0\|_2 \|\hat{x}+x_0\|_2 \leq r_2^*:=\max\{r_0(c_1),r_2(c_2/\sigma \sqrt{\log N})\}.
$$
When $\norm{w}_\infty<\infty$ the term $\sigma\sqrt{\log N}$ may be replaced by $\norm{w}_\infty$.

\vskip0.7cm

The upper estimate of $\max\{r_0,r_2\}$ in Theorem A represents two ranges of noise. It follows from the definition of the fixed points that $r_0$ is dominant if $\sigma \leq r_0/\sqrt{\log N}$. As explained in \cite{Lec-Men} for linear regression, $r_0$ captures the difficulty of recovery in the noise free case, when the only reason for errors is that there are several far-away functions in the class that coincide with the target on the noiseless data. When the noise level $\sigma$ surpasses that threshold, errors occur because of the interaction class members have with the noise, and the dominating term becomes $r_2$.

Of course, there are cases in which $r_0=0$ for $N$ sufficiently large. This is precisely when exact recovery is possible in the noise-free environment. And, in such cases, the error of ERM tends to zero with $\sigma$.

The behavior of ERM in the noise-free case is one of the distinguishing features of sets with well behaved `global complexity' -- because $E$ is not too large. Since $E_R \leq E$ for every $R>0$, it follows that when $N\gtrsim E^2$, $r_0=0$ and that $r_2(\gamma) \leq E/(\gamma\sqrt{N})$. Therefore, on the event from Theorem A,
$$
\|\hat{x}-x_0\|_2 \|\hat{x}+x_0\|_2 \lesssim \sigma \frac{E}{\sqrt{N}} \sqrt{\log N}.
$$
This estimate suffices for many applications. For example, when $T$ is the set of $d$-sparse vectors, one may show (see, e.g. \cite{EM12}) that
$$
E \lesssim \sqrt{d \log(en/d)}.
$$
Hence, by Theorem A, when $N\gtrsim d\log\big(en/d\big)$, with high probability,
$$
\|\hat{x}-x_0\|_2 \|\hat{x}+x_0\|_2 \lesssim \sigma \sqrt{\frac{d \log (en/d)}{N}} \sqrt{\log N}.
$$
The proof of this observation regarding $d$-sparse vectors, and that this estimate is sharp in the minimax sense (up to the logarithmic term) may be found in Section \ref{sec:examples}.

One should note that Theorem~A improves the main result from \cite{EM12} in three ways. First of all, the error rate (the estimate on $\|\hat{x}-x_0\|_2 \|\hat{x}+x_0\|_2$) established in Theorem A is  $\sim E/\sqrt{N}$ (up to logarithmic factors), whereas in \cite{EM12}, it scaled like $c/N^{1/4}$ for very large values of $N$. Second, the error rate scales linearly in the noise level $\sigma$ in Theorem~A. On the other hand, the rate obtained in \cite{EM12} does not decay with $\sigma$ for $\sigma \leq 1$. Finally, the probability estimate has been improved, though it is still likely to be suboptimal.

Although the main motivation for \cite{EM12} was dealing with phase retrieval for sparse classes, and for which Theorem A  is well suited, we next turn to the question of more general classes, the most important example of which is a convex, centrally-symmetric class. For such a class, the global localization is simply too big to yield a good bound.

\begin{Definition} \label{def:fixed-point-localized}
Let
$$
r_N^*(Q)=\inf \left\{ r>0: \ell(T \cap r B_2^n) \leq Q r \sqrt{N}\right\},
$$
$$
s_N^*(\eta)=\inf \left\{s>0: \ell(T \cap s B_2^n) \leq \eta s^2 \sqrt{N}\right\},
$$
and
$$
v_N^*(\zeta)=\inf \left\{v>0: \ell(T \cap v B_2^n) \leq \zeta v^3 \sqrt{N}\right\},
$$
\end{Definition}

The parameters $r_N^*$ and $s_N^*$ have been used in \cite{Lec-Men} to obtain a sharp estimate on the performance of ERM for linear regression in an arbitrary convex set, and relative to $L$-subgaussian measurements. This result is formulated below in a restricted context, analogous to the phase retrieval setup: a linear model $z=\inr{a,x_0}+w$, for an isotropic, $L$-subgaussian vector $a$, independent noise $w$ and $x_0 \in T$.

Let $\hat{x}$ be the output of ERM using the data $(a_i,z_i)_{i=1}^N$ and set $\|w\|_{\psi_2}=\sigma$.

\begin{Theorem} \label{thm:LM-linear}
For every $L \geq 1$ there exist constants $c_1,c_2,c_3$ and $c_4$ that depend only on  $L$ for which the following holds. Let $T\subset\R^d$ be a convex set, put $\eta= c_1/\sigma$ and set $Q = c_2$.
\begin{description}
\item{1.} If $\sigma \geq c_3r_N^*(Q)$ then with probability at least $1-4\exp(-c_4N\eta^2 (s_N^*(\eta))^2)$,
$$
\|x-x_0\|_2 \leq  s_N^*(\eta).
$$
\item{2.} If $\sigma \leq c_3r_N^*(Q)$ then with probability at least $1-4\exp(-c_4NQ^2)$,
$$
\|x-x_0\|_2 \leq r_N^*(Q).
$$
\end{description}
\end{Theorem}
Our main result is a phase retrieval version of Theorem \ref{thm:LM-linear}.
\vskip0.3cm
\noindent{\bf Theorem B.} For every $L \geq 1$, $\kappa_0>0$ and $\beta$ there exist constants $c_1,c_2,c_3$, $c_4,c_5$ and $Q$ that depend only on $L$ and $\kappa_0$ and $\beta$ for which the following holds. Let $T\subset\R^d$ be a convex, centrally-symmetric set, and let $a$ and $w$ be as in Theorem A. Assume that $(\sigma/\|x_0\|_2) \geq c_0r_N^*(Q)/\sqrt{\log N}$, set $\eta=c_1\|x_0\|_2/(\sigma \sqrt{\log N})$ and let $\zeta=c_1/(\sigma\sqrt{\log N})$.
\begin{description}
\item{1.} If $\|x_0\|_2 \geq v_N^*(c_2)$, then with probability at least $1-2\exp(-c_3N\eta^2  (s_N^*(\eta))^2)-2N^{-\beta+1}$,
    $$
    \min\{\|x-x_0\|_2, \|x+x_0\|_2 \} \leq c_4 s_N^*(\eta).
    $$
\item{2.} If $\|x_0\|_2 \leq v_N^*(c_2)$ then with probability at least $1-2\exp(-c_3N\zeta^2  (v_N^*(\zeta))^2)-2N^{-\beta+1}$,
$$
    \max\{\|x-x_0\|_2, \|x+x_0\|_2 \} \leq c_4 v_N^*(\zeta).
    $$
\end{description}
If $(\sigma/\|x_0\|_2) \geq c_0r_N^*(Q)/\sqrt{\log N}$ the same assertion as in 1. and 2. holds, with an upper bound of $r_N^*(Q)$ replacing $s_N^*(\eta)$ and $v_N^*(\zeta)$.

\vskip0.7cm

Theorem B follows from Theorem A and a more transparent description of the localized sets $T_{-,R}(x_0)$ and $T_{+,R}(x_0)$ (see Lemma \ref{localized-sets-large-norm}).

\vskip0.3cm

To put Theorem B in some perspective, observe that $v_N^*$ tends to zero. Indeed, since $\ell(T \cap r B_2^n) \leq \ell(T)$, it follows that $v_N^*(\zeta) \leq (\ell(T)/\sqrt{N} \zeta)^{1/3}$. Hence, for the choice of $\zeta\sim (\sigma \sqrt{\log N})^{-1}$ as in Theorem B,
$$
v_N^* \leq \left(\sigma \ell(T) \sqrt{\frac{\log N}{N}}\right)^{1/3},
$$
which tends to zero when $\sigma \to 0$ and when $N \to \infty$. Therefore, if $x_0 \not = 0$, the first part of Theorem~B describes the `long term' behaviour of ERM.

Also, and using the same argument,
$$
r_N^*(Q) \leq \frac{\ell(T)}{Q \sqrt{N}}.
$$
Thus, for every $\sigma>0$ the problem becomes `high noise' in the sense that the condition $(\sigma/\|x_0\|_2) \geq c_0r_N(Q)/\sqrt{\log N}$ is satisfied when $N$ is large enough.

In the typical situation, which is both `high noise' and `large $\|x_0\|_2$', the error rate depends on $\eta=c_1\|x_0\|_2/\sigma \sqrt{\log N}$. We believe that the $1/\sqrt{\log N}$ factor is an artifact of the proof, but the other term, $\|x_0\|_2/\sigma$ is the signal-to-noise ratio, and is rather natural.


\vskip0.3cm

Although Theorem A and Theorem B clearly improve the results from \cite{EM12}, it is natural to ask whether these are optimal in a more general sense. The final result presented here is that Theorem~B is close to being optimal in the minimax sense.  The formulation and proof of the minimax lower bound is presented in Section \ref{sec:minimax}.

Finally, we end the article with two examples of classes that are of interest in phase retrieval: $d$-sparse vectors and the unit ball in $\ell_1^n$. The first is a class with a fixed `local complexity', and the second has a growing `local complexity'.

\section{Preliminaries} \label{sec:pre}
Throughout this article, absolute constants are denoted by $C,c,c_1,...$ etc. Their value may change from line to line. The fact that there are absolute constants $c,C$ for which $ca \leq b \leq Ca$ is denoted by $a \sim b$; $a \lesssim b$ means that $a \leq cb$, while $a \sim_L b$ means that the constants depend only on the parameter $L$.

For $1 \leq p \leq \infty$, let $\|\cdot\|_p$ be the $\ell_p^n$ norm endowed on $\R^n$, and for a function $f$ (or a random variable $X$) on a probability space, set $\|f\|_{L_p}$ to be its $L_p$ norm.

Other norms that play a significant role here are the Orlicz norms. For basic facts on these norms we refer the reader to \cite{LT,VW}.

Recall that for $\alpha \geq 1$,
$$
\|f\|_{\psi_\alpha} = \inf \{c>0 : \E \exp(|f|^\alpha/c^\alpha) \leq 2\},
$$
and it is straightforward to extend the definition for $0<\alpha<1$.

Orlicz norms measure the rate of decay of a function. One may verify that $\|f\|_{\psi_\alpha} \sim \sup_{p \geq 1} \|f\|_{L_p}/p^{1/\alpha}$. Moreover, for $t \geq 1$, $Pr(|f| \geq t) \leq 2 \exp(-ct^\alpha/ \|f\|_{\psi_\alpha}^\alpha)$, and $\|f\|_{\psi_\alpha}$ is equivalent to the smallest constant $\kappa$ for which $Pr(|f| \geq t) \leq 2 \exp(-t^\alpha/ \kappa^\alpha)$ for every $t \geq 1$.

\begin{Definition} \label{def-subgauss-pre}
A random variable is $L$-subgaussian if it has a bounded $\psi_2$ norm and $\|X\|_{\psi_2} \leq L \|X\|_{L_2}$.
\end{Definition}
Observe that for $L$-subgaussian random variables, all the $L_p$ norms are equivalent and their tails exhibits a faster decay than the corresponding gaussian. Indeed, if $X$ is $L$-subgaussian,
$$
\|X\|_{L_p} \lesssim \sqrt{p}\|X\|_{\psi_2} \lesssim L\sqrt{p}\|X\|_{L_2},
$$
and for every $t \geq 1$,
$$
Pr (|X| > t) \leq 2\exp(-ct^2/\|X\|_{\psi_2}^2)\leq 2\exp(-ct^2/(L^2\|X\|_{L_2}^2))
$$
for a suitable absolute constant $c$.

It is standard to verify that for every $f,g$, $\|f g\|_{\psi_1} \lesssim \|f\|_{\psi_2} \|g\|_{\psi_2}$, and that if $X_1,...,X_N$ are independent copies of $X$ and $1 \leq \alpha \leq 2$, then
\begin{equation} \label{eq:Pisier-ineq}
\|\max_{1 \leq i \leq N} X_i \|_{\psi_\alpha} \lesssim  \|X\|_{\psi_\alpha} \log^{1/\alpha} N.
\end{equation}

An additional feature of $\psi_\alpha$ random variables is concentration, namely that if $(X_i)_{i=1}^N$ are independent copies of a $\psi_\alpha$ random variable $X$, then $N^{-1}\sum_{i=1}^N X_i$ concentrates around $\E X$. One example of such a concentration result is the following Bernstein-type inequality (see, e.g., \cite{VW}).
\begin{Theorem} \label{thm:Bernstein}
There exists an  absolute constant $c_0$ for which the following holds. If $X_1,...,X_N$ are independent copies of a $\psi_1$ random variable $X$, then for every $t>0$,
$$
Pr\left(\left|\frac{1}{N}\sum_{i=1}^N X_i -\E X\right| > t\|X\|_{\psi_1} \right) \leq
2\exp(-c_0N \min\{t^2,t\}).
$$
\end{Theorem}
One important example of a probability space considered here is the discrete space $\Omega=\{1,...,N\}$, endowed with the uniform probability measure. Functions on $\Omega$ can be viewed as vectors in $\R^N$ and the corresponding $L_p$ and $\psi_\alpha$ norms are denoted by $\| \cdot \|_{L_p^N}$ and $\| \cdot \|_{\psi_\alpha^N}$.

A significant part of the proofs of Theorem A  has to do with the behaviour of a monotone non-increasing rearrangement of vectors. Given $v \in \R^N$, let $(v_i^*)_{i=1}^N$ be a non-increasing rearrangement of $(|v_i|)_{i=1}^N$. It turns out that the $\psi_\alpha^N$ norm captures  information on the coordinates of $(v_i^*)_{i=1}^N$.

\begin{Lemma} \label{lemma-psi-2-discrete}
For every $1 \leq \alpha \leq 2$ there exist constants $c_1$ and $c_2$ that depend only on $\alpha$ for which the following holds. For every $v \in \R^N$,
$$
c_1\sup_{i \leq N} \frac{v_i^*}{\log^{1/\alpha}(eN/i)} \leq \|v\|_{\psi_\alpha^N} \leq c_2\sup_{i \leq N} \frac{v_i^*}{\log^{1/\alpha}(eN/i)}.
$$
\end{Lemma}
\proof
We will prove the claim only for $\alpha=2$ as the other cases follow an identical path.

Let $v \in \R^N$ and denote by $Pr$ the uniform probability measure on $\Omega=\{1,\ldots,N\}$. By the tail characterization of the $\psi_2$ norm,
$$
N^{-1} |\{j: |v_j| >t\}| = Pr (|v| >t ) \leq 2\exp(-ct^2/\|v\|_{\psi_2^N}^2).
$$
Hence, for $t_i=c^{-1/2}\|v\|_{\psi_2^N} \sqrt{\log(eN/i)}$, $|\{j: |v_j| >t_i\}| \leq 2i/e \leq i$, and for every $1 \leq i \leq N$, $v_i^* \leq t_i$. Therefore,
$$
\sup_{i \leq N} \frac{v_i^*}{\sqrt{\log(eN/i)}} \leq c^{-1/2}\|v\|_{\psi_2^N},
$$
as claimed.

In the reverse direction, consider
$$
{\cal B}=\big\{\beta>0 : \forall \ 1 \leq i \leq N, \ \|v\|_{\psi_2^N} \geq \beta v_i^*/\sqrt{\log(eN/i)} \big\}.
$$
It is enough to show that ${\cal B}$ is bounded by a constant that is independent of $v$. To that end, fix $\beta \in {\cal B}$ and without loss of generality, assume that $\beta>2$. Set $B=\sup_{i \leq N} \beta v_i^*/\sqrt{\log(eN/i)}$ and since $\beta \in {\cal B}$, $\|v\|_{\psi_2^N} \geq B$.

Also, since $1/\beta^2 <1$,
$$
\sum_{i=1}^N \left(\frac{1}{i}\right)^{1/\beta^2} \leq 1+\int_1^N \left(\frac{1}{x}\right)^{1/\beta^2}dx \leq \frac{N^{1-1/\beta^2}}{1-1/\beta^2}.
$$
Therefore,
%
\begin{align*}
&\sum_{i=1}^N \exp(v_i^2/B^2) =  \sum_{i=1}^N \exp((v_i^*)^2/B^2) \leq \sum_{i=1}^N \exp(\beta^{-2}\log(eN/i))
\\
\leq & \sum_{i=1}^N \left(\frac{eN}{i}\right)^{1/\beta^2} \leq (eN)^{1/\beta^2} \cdot \frac{N^{1-1/\beta^2}}{1-1/\beta^2}
\leq   \frac{N e^{1/\beta^2}}{1-1/\beta^2} < 2N,
\end{align*}
provided that $\beta \geq c_1$. Thus, if $\beta \geq c_1$, $\|v\|_{\psi_2^N}<B$ which is a contradiction, showing that ${\cal B}$ is bounded by $c_1$.
\endproof

\subsection{Empirical and Subgaussian processes} \label{sec:emp-subgauss}
The sampling method used here is isotropic and $L$-subgaussian, meaning that the vectors
$(a_i)_{i=1}^N$ are independent and distributed according to a probability measure $\mu$ on $\R^n$ that is both isotropic and $L$-subgaussian \cite{VW}:
\begin{Definition} \label{def-iso-subgaussian}
Let $\mu$ be a probability measure on $\R^n$ and let $a$ be distributed according to $\mu$. The measure $\mu$ is isotropic if for every $t \in \R^n$, $\E\inr{a,t}^2 = \|t\|_{2}^2$. It is  $L$-subgaussian if for every $t \in \R^n$ and every $u \geq 1$, $Pr(|\inr{a,t}| \geq Lu\|\inr{t,a}\|_{2}) \leq 2\exp(-u^2/2)$.
\end{Definition}

Given $T \subset \R^n$, let $d_T = \sup_{t \in T} \|t\|_2$ and put $k_*(T)= (\ell(T)/d_T)^2$. The latter appears naturally in the context of Dvoretzky type theorems, and in particular, in Milman's proof of Dvoretzky's Theorem (see, e.g., \cite{MS}).

\begin{Theorem} \label{thm:monotone-class} \cite{Men-subgauss}
For every $L \geq 1$ there exist constants $c_1$ and $c_2$ that depend only on $L$ and for which the following holds. For every $u \geq c_1$, with probability at least $1-2\exp(-c_2u^2k_*(T))$, for every $t \in T$ and every $I \subset \{1,...,N\}$,
$$
\left(\sum_{i \in I} \inr{t,a_i}^2 \right)^{1/2} \leq Lu^3 \left(\ell(T) + d_T\sqrt{|I|\log(eN/|I|)}\right).
$$
\end{Theorem}
For every integer $N$, let $j_{T}$ be the largest integer $j$ in $\{1,...,N\}$ for which
$$
\ell(T) \geq d_T \sqrt{j \log (eN/j)}.
$$
It follows from Theorem \ref{thm:monotone-class} that if $t \in T$ and $|I| \leq j_{T}$,
$$
(\sum_{i \in I} \inr{t,a_i}^2)^{1/2} \lesssim_{L,u} \ell(T),
$$
and if $|I| \geq j_{T}$,
$$
(\sum_{i \in I} \inr{t,a_i}^2)^{1/2} \lesssim_{L,u} d_T\sqrt{|I|\log(eN/|I|)}.
$$
 Therefore, if $v=(\inr{t,a_i})_{i=1}^N$ and $(v_i^*)_{i=1}^N$ is a monotone non-increasing rearrangement of $(|v_i|)_{i=1}^N$, then
\begin{equation} \label{eq:monotone-decomp}
v_i^* \leq \left(\frac{1}{i}\sum_{j=1}^i (v_j^*)^2 \right)^{1/2} \lesssim_{L,u}
\left\{\begin{array}{cc}
\frac{\ell(T)}{\sqrt{i}} & \mbox{if} \ \ i \leq j_T\\
&\\
d_T \sqrt{\log(eN/i)} & \mbox{otherwise}.
\end{array}\right.
\end{equation}
This observation will be used extensively in what follows.

The next fact deals with product processes.
\begin{Theorem} \label{thm:emp-subgaussian-prod} \cite{Men-subgauss}
There exist absolute constants $c_0, c_1$ and $c_2$ for which the following holds. If $T_1,T_2 \subset \R^n$, $1 \leq 2^j \leq N$ and $u \geq c_0$, then with probability at least $1-2\exp(-c_1 u^2 2^j)$,
\begin{align*}
& \sup_{t \in T_1, \ s \in T_2} \left|\sum_{i=1}^N \inr{a_i,t}\inr{a_i,s} - \E \inr{a,t}\inr{a,s}\right|
\\
\leq & c_2L^2 u^2 \left(\ell(T_1) \ell(T_2)+u \sqrt{N} \left(\ell(T_1)d(T_2)+\ell(T_2)d(T_1)+2^{j/2}d(T_1)d(T_2)\right) \right)
\end{align*}
and
\begin{align*}
& \sup_{t \in T_1, \ s \in T_2} \left|\sum_{i=1}^N |\inr{a_i,t}\inr{a_i,s}| - \E |\inr{a,t}\inr{a,s}|\right|
\\
\leq & c_2L^2 u^2 \left(\ell(T_1) \ell(T_2)+u \sqrt{N} \left(\ell(T_1)d(T_2)+\ell(T_2)d(T_1)+2^{j/2}d(T_1)d(T_2)\right) \right).
\end{align*}
\end{Theorem}
\begin{Remark} \label{rem:Bernoulli-product}
Let $(\eps_i)_{i=1}^N$ be independent, symmetric, $\{-1,1\}$-valued random variables. It follows from the results in \cite{Men-subgauss} that with the same probability estimate in Theorem \ref{thm:emp-subgaussian-prod} and relative to the product measure $(\eps \otimes X)^N$,
\begin{align*}
& \sup_{t \in T_1, \ s \in T_2} \left|\sum_{i=1}^N \eps_i\inr{a_i,t}\inr{a_i,s}\right|
\\
\lesssim & L^2 u^2 \left(\ell(T_1) \ell(T_2)+u \sqrt{N} \left(\ell(T_1)d(T_2)+\ell(T_2)d(T_1)+2^{j/2}d(T_1)d(T_2)\right) \right).
\end{align*}
\end{Remark}
\vskip0.5cm

Assume that $(k^*(T_1))^{1/2}=\ell(T_1)/d(T_1) \geq \ell(T_2)/d(T_2)$. Setting $2^{j/2}=\ell(T_1)/d(T_1)$, Theorem \ref{thm:emp-subgaussian-prod} and Remark \ref{rem:Bernoulli-product} yield that with probability at least $1-2\exp(-c_1u^2k_*(T_1))$,
\begin{equation*}
\sup_{t \in T_1, \ s \in T_2} \left|\sum_{i=1}^N \inr{a_i,t}\inr{a_i,s} - \E \inr{a,t}\inr{a,s}\right|
\lesssim L^2 u^2 \ell(T_1)\left( \ell(T_2)+u \sqrt{N}d(T_2)\right),
\end{equation*}
\begin{equation*}
\sup_{t \in T_1, \ s \in T_2} \left|\sum_{i=1}^N |\inr{a_i,t}\inr{a_i,s}| - \E |\inr{a,t}\inr{a,s}|\right|
\lesssim L^2 u^2 \ell(T_1)\left( \ell(T_2)+u \sqrt{N}d(T_2)\right)
\end{equation*}
and
\begin{equation} \label{eq:Bernoulli-prod-est-useful}
\sup_{t \in T_1, \ s \in T_2} \left|\sum_{i=1}^N \eps_i\inr{a_i,t}\inr{a_i,s}\right|
\lesssim L^2 u^2 \ell(T_1)\left( \ell(T_2)+u \sqrt{N}d(T_2)\right).
\end{equation}

One case which is of particular interest is when $T_1=T_2=T$, and then,
with probability at least $1-2\exp(-c_1u^2k_*(T))$,
\begin{equation*}
\sup_{t \in T} \left|\sum_{i=1}^N \inr{a_i,t}^2 - \E \inr{a,t}^2\right|
\lesssim L^2 u^2 \left(\ell^2(T) +u \sqrt{N}d(T)\ell(T)\right).
\end{equation*}

\subsection{Monotone rearrangement of coordinates} \label{sec:monotone}
The first goal of this section is to investigate the coordinate structure of $v \in \R^m$, given information on its norm in various $L_p^m$ and $\psi_\alpha^m$ spaces. The vectors we will be interested in are of the form $(\inr{a_i,t})_{i=1}^N$ for $t \in T$, and for which, thanks to the results presented in Section \ref{sec:emp-subgauss}, one has the necessary information at hand.

It is standard to verify that if $\|v\|_{\psi_\alpha^m} \leq A$, then $\|v\|_p \lesssim_{p,\alpha}  A \cdot m^{1/p}$. Thus, $\|v\|_{L_p^m} \lesssim_p \|v\|_{\psi_\alpha^m}$.

It turns out that if the two norms are equivalent, $v$ is regular in some sense. The next lemma, which is a version of the Paley-Zygmund Inequality, (see, e.g. \cite{GdlP}), describes the regularity properties needed here in the case $p=\alpha=1$.

\begin{Lemma} \label{lemma:comp-emp-norms}
For every $\beta>1$ there exist constants $c_1$ and $c_2$ that depend only on $\beta$ and for which the following holds. If $\|v\|_{\psi_1^m} \leq \beta \|v\|_{L_1^m}$, there exists $I \subset \{1,...,m\}$ of cardinality at least $c_1 m$, and for every $i \in I$, $|v_i| \geq c_2 \|v\|_{L_1^m}$.
\end{Lemma}

\proof
Recall that $\|v\|_{\psi_1^m} \sim \sup_{1 \leq i \leq m} v_i^*/\log(em/i)$. Hence,
for every $1 \leq j \leq m$,
$$
\sum_{\ell =1}^j v_\ell^* \lesssim \|v\|_{\psi_1^m} \sum_{\ell =1}^j \log(em/\ell) \lesssim \beta \|v\|_{L_1^m} j\log(em/j).
$$
Therefore,
\begin{equation*}
m\|v\|_{L_1^m}  =    \sum_{\ell=1}^m |v_\ell| = \sum_{\ell \leq j} v_\ell^* + \sum_{\ell=j+1}^m v_\ell^*
\leq  c_0\beta \|v\|_{L_1^m} j\log(em/j) + \sum_{\ell=j+1}^m v_\ell^*.
\end{equation*}
Setting $c_1(\beta) \sim 1/(\beta \log(e\beta))$ and $j=c_1(\beta)m$,
$$
c_0 \beta \|v\|_{L_1^m}j\log(em/j) \leq (m/2)\|v\|_{L_1^m}.
$$
Thus,
$\sum_{\ell=j+1}^m v_\ell^* \geq (m/2) \|v\|_{L_1^m}$, while
$$
v_{j+1}^* \leq \frac{1}{j+1}\sum_{\ell \leq j+1} v_\ell^* \lesssim \beta \log (e\beta) \|v\|_{L_1^m}.
$$

Let $I$ be the set of the $m-j$ smallest coordinates of $v$. Fix $\eta>0$ to be named later, put $I_\eta \subset I$ to be the set of coordinates in $I$ for which $|v_i| \geq \eta \|v\|_{L_1^m}$ and denote by $I_\eta^c$ its complement in $I$. Therefore,
\begin{align*}
(m/2) \|v\|_{L_1^m} \leq & \sum_{\ell \geq j+1} v_\ell^* = \sum_{\ell \in I_\eta} |v_\ell| + \sum_{\ell \in I_\eta^c} |v_\ell| \leq v_{j+1}^* |I_\eta| + \eta \|v\|_{L_1^m} |I_\eta^c|
\\
\lesssim & \|v\|_{L_1^m} |I| \left( \beta \log (e\beta)  \frac{|I_\eta|}{|I|} + \eta\frac{|I_\eta^c|}{|I|}\right).
\end{align*}
Hence,
\begin{equation*}
\frac{m}{2} \lesssim  |I| \left(\beta \log(e \beta) \frac{|I_\eta|}{|I|} + \eta \left(1-\frac{|I_\eta|}{|I|}\right) \right)
\lesssim  m \left( \left(\beta \log(e \beta)-\eta\right)\frac{|I_\eta|}{|I|} + \eta\right).
\end{equation*}
If $\eta=\min\{1/4, (\beta/2)\log(e \beta)\}$, then $|I_\eta| \geq (\eta/2)|I| \geq c_2(\beta)m$, as claimed.
\endproof

Next, let us turn to decomposition results for vectors of the form $(\inr{a_i,t})_{i=1}^N$. Recall that for a set $T \subset \R^N$, $j_T$ is the largest integer for which $\ell(T) \geq d_T \sqrt{j \log (eN/j)}$.
\begin{Lemma} \label{lemma:decomposition-ratio-set}
For every $L>1$ there exist constants $c_1$ and $c_2$ that depend only on $L$ and for which the following holds.
Let $T \subset \R^n$ and set $W=\{t/\|t\|_2 : t \in T\} \subset S^{n-1}$. With probability at least $1-2\exp(-c_1\ell^2(W))$, for every $t \in T$, $(\inr{a_i,t})_{i=1}^N=v_1+v_2$ and $v_1,v_2$ have the following properties:
\begin{description}
\item{1.} The supports of $v_1$ and $v_2$ are disjoint.
\item{2.} $\|v_1\|_2 \leq c_2\ell(W)\|t\|_2$ and $|{\rm supp}(v_1)| \leq j_{W}$.
\item{3.} $\|v_2\|_{\psi_2^N} \leq c_2 \|t\|_2$.
\end{description}
\end{Lemma}
\proof
Fix $t \in T$ and let $J_t \subset \{1,...,N\}$ be the set of the largest $j_{W}$ coordinates of $(|\inr{a_i,t}|)_{i=1}^N$. Set
$$
{\bar v}_1=(\inr{a_j,t/\|t\|_2})_{j \in J_t} \ \ {\rm and} \ \ {\bar v}_2=(\inr{a_j,t/\|t\|_2})_{j \in J_t^c}.
$$
By Theorem \ref{thm:monotone-class} and the characterization of the $\psi_2^N$ norm of a vector using the monotone rearrangement of its coordinates (Lemma \ref{lemma-psi-2-discrete}),
$$
\|\bar{v}_1\|_2 \lesssim L\ell(W),  \ \ {\rm and} \ \ \|\bar{v}_2\|_{\psi_2^N} \lesssim L.
$$
To conclude the proof, set $v_1=\|t\|_2 \bar{v}_1$ and $v_2=\|t\|_2 {\bar v}_2$.
\endproof

Recall that for every $R>0$,
$$
T_{+,R} =\left\{ \frac{t+s}{\|t+s\|_2} : \ t,s \in T, \ \|t+s\|_2 \|t-s\|_2 \geq R \right\},
$$
and a similar definition holds for $T_{-,R}$. Set $j_{+,R}=j_{T_{+,R}}$, $j_{-,R}=j_{T_{+,R}}$ and $E_R = \max\{ \ell(T_{+,R}), \ell(T_{-,R}) \}$. Combining the above estimates leads to the following corollary.
\begin{Corollary} \label{cor:all-info-on-proj}
For every $L>1$ there exist constants $c_1,c_2,c_3$ and $c_4$ that depend only on $L$ for which the following holds. Let $T \subset \R^n$ and $R>0$, and consider $T_{+,R}$ and $T_{-,R}$ as above. With probability at least $1-4\exp(-c_1L^2\min\{\ell^2(T_{+,R}),\ell^2(T_{-,R})\})$, for every $s,t \in T$ for which $\|t-s\|_2 \|t+s\|_2 \geq R$,
\begin{description}
\item{1.} $(\inr{s-t,a_i})_{i=1}^N = v_1+v_2$, for vectors $v_1$ and $v_2$ of disjoint supports satisfying
$$
|{\rm supp}(v_1)| \leq j_{-,R}, \ \ \|v_1\|_2 \leq c_2 \ell(T_{-,R}) \|s-t\|_2 \ \ {\rm and} \ \ \|v_2\|_{\psi_2^N} \leq c_2 \|s-t\|_2.
$$
\item{2.} $(\inr{s+t,a_i})_{i=1}^N = u_1+u_2$, for vectors $u_1$ and $u_2$ of disjoint supports satisfying
$$
|{\rm supp}(u_1)| \leq j_{+,R}, \ \ \|u_1\|_2 \leq c_2 \ell(T_{+,R}) \|s+t\|_2 \ \ {\rm and} \ \ \|u_2\|_{\psi_2^N} \leq c_2 \|s+t\|_2.
$$
\item{3.} If $h_{s,t}(a)=\inr{\frac{s+t}{\|s+t\|_2},a}\inr{\frac{s-t}{\|s-t\|_2},a}$, then
\begin{equation*}
 \left|\frac{1}{N}\sum_{i=1}^N |h_{s,t}(a_i)| - \E |h_{s,t}| \right|
\leq c_3 \left(\frac{E_R}{\sqrt{N}}+\frac{E_R^2}{N}\right).
\end{equation*}
\end{description}
In particular, recalling that for every $s,t \in T$,
\begin{equation*}
\E |\inr{s+t,a}\inr{s-t,a}| \geq \kappa_0 {\|s+t\|_2\|s-t\|_2},
\end{equation*}
it follows that if $\sqrt{N} \geq c_4(L)E_R/\kappa_0$ then
\begin{description}
\item{4.}
\begin{equation} \label{eq:property-3-prime}
\frac{\kappa_0}{2}\|s+t\|_2\|s-t\|_2 \leq \frac{1}{N}\sum_{i=1}^N |\inr{s+t,a_i}\inr{s-t,a_i}| \lesssim_L  \|s+t\|_2\|s-t\|_2.
\end{equation}
\end{description}
\end{Corollary}

From here on, denote by $\Omega_{1,R}$ the event on which Corollary \ref{cor:all-info-on-proj} holds for the sets $T_{+,R}$ and $T_{-,R}$ and samples of cardinality $N \gtrsim_{L} E_R^2/\kappa_0^2$.

\begin{Lemma} \label{lemma:large-coordinates}
There exist constants $c_0$ depending only on $L$ and $c_1, \kappa_1$ that depend only on $\kappa_0$ and $L$ for which the following holds. If $N \geq c_0 E_R^2/\kappa_0^2$, then for $(a_i)_{i=1}^N \in \Omega_{1,R}$, for every $s,t \in T$ for which $\|s-t\|_2 \|s+t\|_2 \geq R$, there is $I_{s,t} \subset \{1,...,N\}$ of cardinality at least $\kappa_1N$, and for every $i \in I_{s,t}$,
$$
|\inr{s-t,a_i} \inr{s+t,a_i}| \geq c_1 \|s-t\|_2 \|s+t\|_2.
$$
\end{Lemma}

Lemma \ref{lemma:large-coordinates} is an empirical ``small-ball" estimate, as it shows that with high probability, and for every pair $s,t$ as above, many of the coordinates of $(|\inr{a_i,s-t}| \cdot |\inr{a_i,s+t}|)_{i=1}^N$ are large.

\proof
Fix $s,t \in T$ as above and set
$$
y=(\inr{s-t,a_i})_{i=1}^N, \ \ {\rm and} \ \ x=(\inr{s+t,a_i})_{i=1}^N.
$$
 Let $y=v_1+v_2$ and $x=u_1+u_2$ as in Corollary \ref{cor:all-info-on-proj}. Let $j_0=\max\{j_{-,R},j_{+,R}\}$ and put $J={\rm supp}(v_1) \cup {\rm supp}(u_1)$. Observe that $|J| \leq 2j_0$ and that
\begin{align*}
& \sum_{j \in J} |y(j)| \cdot |x(j)| \leq  \sum_{j \in {\rm supp}(v_1)} |v_1(j) x(j)| + \sum_{j \in {\rm supp}(u_1)} |y(j)u_1(j)|
\\
\leq & \|v_1\|_2 \left(\sum_{i=1}^{2j_0} (x^2(j))^* \right)^{1/2} + \|u_1\|_2 \left(\sum_{i=1}^{2j_0} (y^2(j))^* \right)^{1/2}
\\
\lesssim_L & \ell(T_{-,R}) \|s-t\|_2 \cdot \sqrt{j_0\log(eN/j_0)} \|s+t\|_2
\\
+ & \ell(T_{+,R}) \|s+t\|_2 \cdot \sqrt{j_0\log(eN/j_0)} \|s-t\|_2
\\
\lesssim_L &  E_R^2 \|s-t\|_2 \|s+t\|_2 \leq \frac{\kappa_0 N}{4} \|s-t\|_2 \|s+t\|_2,
\end{align*}
because, by the definition of $j_0$, $\sqrt{j_0\log(eN/j_0)} \lesssim \max\{\ell(T_{-,R}),\ell(T_{+,R})\}$ and since $N \geq c_0 E_R^2/\kappa_0^2$ for $c_0=c_0(L)$ large enough.

Thus, by \eqref{eq:property-3-prime},
$$
\sum_{j \in J^c} |y(j)x(j)| \geq N\kappa_0\|s-t\|_2\|s+t\|_2/4.
$$
Set $m=|J^c|$ and let $z=(y(j)x(j))_{j \in J^c}=(v_2(j)u_2(j))_{j \in J^c}$. Since $N \gtrsim_{L} E_R^2/\kappa_0^2$, it is evident that $j_0 \leq N/2$; thus $N/2 \leq m \leq N$ and
$$
\|z\|_{L_1^m} = \frac{1}{m} \sum_{j \in J^c} |y(j)x(j)| \geq \frac{N}{4m} \kappa_0 \|s-t\|_2\|s+t\|_2 \gtrsim \kappa_0 \|s-t\|_2\|s+t\|_2.
$$
On the other hand,
$$
\|z\|_{\psi_1^m} \leq \|(v_2u_2(j))_{j \in J_c}\|_{\psi_1^m} \lesssim \|v_2\|_{\psi_2^m} \|u_2\|_{\psi_2^m} \lesssim_L  \|s-t\|_2 \|s+t\|_2,
$$
and $z$ satisfies the assumption of Lemma \ref{lemma:comp-emp-norms} for $\beta=c_1(L,\kappa_0)$. The claim follows immediately from that lemma.
\endproof

\section{Proof of Theorem A} \label{sec:squared-loss}
It is well understood that when analyzing properties of ERM relative to a loss $\ell$, studying the excess loss functional is rather natural. The excess loss shares the same empirical minimizer as the loss, but it has additional qualities: for every $x \in T$, $\E {\cal L}_x \geq 0$ and ${\cal L}_{x_0} = 0$.

Since $0$ is a potential minimizer of $\{P_N {\cal L}_x : x \in T\}$, the minimizer $\hat{x}$ satisfies that $P_N {\cal L}_{\hat{x}} \leq 0$, giving one a way of excluding parts of $T$ as potential empirical minimizers. One simply has to show that with high probability, those parts belong to the set $\{x : P_N {\cal L}_x >0\}$, for example, by showing that $P_N {\cal L}_x$ is equivalent to $\E {\cal L}_x$, as the latter is positive for points that are not true minimizers.

The squared excess loss has a simple decomposition to two processes: a quadratic process and a multiplier one. Indeed, given a class of functions $F$ and $f \in F$,
$$
\big(f(a)-y\big)^2-\big(f^*(a)-y\big)^2=\big(f(a)-f^*(a)\big)^2-2\big(f(a)-f^*(a)\big)\big(f^*(a)-y\big).
$$
where, as always, $f^*$ is a minimizer of the functional $\E\big(f(a)-y\big)^2$ in $F$.

In the phase retrieval problem, $y=\inr{x_0,a}^2+w$ for a noise $w$ that is independent of $a$, and each $f_x \in F$ is given by $f_x=\inr{x,\cdot}^2$. Thus,
\begin{align*}
\cL_x(a,y)&=\ell_x(a,y)-\ell_{x_0}(a,y)=\big(f_x(a)-y\big)^2-\big(f_{x_0}(a)-y\big)^2\\
&=\left(\inr{x-x_0,a} \inr{x+x_0,a}\right)^2-2w\inr{x-x_0,a} \inr{x+x_0,a}.
\end{align*}

Since $w$ is a mean-zero random variable that is independent of $a$, and by Assumption \ref{ass:small-ball},
$$
\E \cL_x(a,y) = \E |\inr{x-x_0,a} \inr{x+x_0,a}| \geq \kappa_0^2 \|x-x_0\|_2^2 \|x+x_0\|_2^2.
$$
Therefore, $\E\big(f_x(a)-y\big)^2$ has a unique minimizer in $F$: $f^*=f_{x_0}=f_{-x_0}$.

To show that $P_N {\cal L}_x>0$ on a large subset $T^\prime \subset T$, it suffices to obtain a high probability lower bound on
$$
\inf_{x \in T^\prime} \frac{1}{N} \sum_{i=1}^N \left(\inr{x-x_0,a_i} \inr{x+x_0,a_i}\right)^2
$$
that dominates a high probability upper bound on
$$
\sup_{x \in T^\prime}  \left| \frac{2}{N} \sum_{i=1}^N w_i  \inr{x-x_0,a_i} \inr{x+x_0,a_i} \right|.
$$
The set $T^\prime$ that will be used is $T_R=\{x \in T : \|x-x_0\|_2 \|x+x_0\|_2 \geq R\}$ for a well-chosen $R$.

\begin{Theorem} \label{thm:quadratic-small-ball}
There exist a constant $c_0$ depending only on $L$, and constants $c_1, \kappa_1$ depending only on $\kappa_0$ and $L$ for which the following holds. For every $R>0$ and $N \geq c_0 E_R^2/\kappa_0^2$, with probability at least $$1-4\exp(-c_1L^2\min\{\ell^2(T_{+,R}),\ell^2(T_{-,R})\}),$$ for every $x \in T_R$,
\begin{equation*}
\frac{1}{N}\sum_{i=1}^N\inr{x_0-x,a_i}^2 \inr{x_0+x,a_i}^2 \geq c_1 \|x_0-x\|_2^2 \|x_0+x\|_2^2.
\end{equation*}
\end{Theorem}
Theorem \ref{thm:quadratic-small-ball} is an immediate outcome of Lemma \ref{lemma:large-coordinates}

\begin{Theorem} \label{thm:multi}
There exist absolute constants $c_1$ and $c_2$ for which the following holds. For every $\beta>1$, with probability at least
$$
1-2\exp(-c_1L^2\min\{\ell^2(T_{+,R}),\ell^2(T_{-,R})\})-2N^{-(\beta-1)},
$$
for every $x \in T_R$,
$$
\left|\frac{1}{N} \sum_{i=1}^N w_i \inr{x-x_0,a_i}\inr{x+x_0,a_i} \right| \leq c_2 \sqrt{\beta}  \|w\|_{\psi_2} \sqrt{\log{N}} \cdot \frac{E_R}{\sqrt{N}} \|x-x_0\|_2 \|x+x_0\|_2.
$$
\end{Theorem}

\proof
By standard properties of empirical process, and since $w$ is mean-zero and independent of $a$, it suffices to estimate
$$
\sup_{x \in T_R} \left|\frac{1}{N} \sum_{i=1}^N \eps_i |w_i| \inr{x-x_0,a_i}\inr{x+x_0,a_i} \right|,
$$
for independent signs $(\eps_i)_{i=1}^N$. By the contraction principle for Bernoulli processes (see, e.g., \cite{LT}), it follows that for every fixed $(w_i)_{i=1}^N$ and $(a_i)_{i=1}^N$,
\begin{align*}
& Pr_\eps \left( \sup_{x \in T_R} \left|\frac{1}{N} \sum_{i=1}^N \eps_i |w_i| \inr{\frac{x-x_0}{\|x-x_0\|_2},a_i}\inr{\frac{x+x_0}{\|x+x_0\|_2},a_i} \right| > u \right)
\\
\leq & 2 Pr_\eps \left( \max_{i \leq N} |w_i| \cdot \sup_{x \in T_R} \left|\frac{1}{N} \sum_{i=1}^N \eps_i \inr{\frac{x-x_0}{\|x-x_0\|_2},a_i}\inr{\frac{x+x_0}{\|x+x_0\|_2},a_i} \right| > \frac{u}{2} \right).
\end{align*}
Applying Remark \ref{rem:Bernoulli-product}, if $N \gtrsim_L E_R$ then with $(\eps \otimes a)^N$-probability of at least $1-2\exp(-c_1L^2\min\{\ell^2(T_{+,R}),\ell^2(T_{-,R})\})$,
$$
 \sup_{x \in T_R} \left|\frac{1}{N} \sum_{i=1}^N \eps_i \inr{\frac{x-x_0}{\|x-x_0\|_2},a_i}\inr{\frac{x+x_0}{\|x+x_0\|_2},a_i} \right| \leq c_2L^2 \frac{E_R}{\sqrt{N}}.
$$
Also, because $w$ is a $\psi_2$ random variable,
$$
Pr(w_1^* \geq t\|w\|_{\psi_2}) \leq 2N \exp(-t^2/2),
$$
and thus, $w_1^* \leq \sqrt{2 \beta \log N} \|w\|_{\psi_2}$ with probability at least $1-2N^{-\beta+1}$.

Combining the two estimates and a Fubini argument, it follows that with probability at least $1-2\exp(-c_1L^2\min\{\ell^2(T_{+,R}),\ell^2(T_{-,R})\})-2N^{-\beta+1}$, for every $x \in T_R$,
$$
\left|\frac{1}{N} \sum_{i=1}^N w_i \inr{x-x_0,a_i}\inr{x+x_0,a} \right| \leq c_3L^2 \sqrt{\beta} \|w\|_{\psi_2} \sqrt{\log{N}} \frac{E_R}{\sqrt{N}} \cdot \|x-x_0\|_2 \|x+x_0\|_2.
$$
\endproof
On the intersection of the two events appearing in Theorem \ref{thm:quadratic-small-ball} and Theorem \ref{thm:multi}, if $N \gtrsim_{\kappa_0,L} E_R^2$ and setting 
$$
\rho=\|x-x_0\|_2 \|x+x_0\|_2 \geq R\geq r_2(c_1\kappa_0/(c_2 L^2 \sqrt{\beta})),
$$ 
then for every $x \in T_R$,
\begin{align*}
P_N {\cal L}_x \geq & \left(c_1 \kappa_0^2 \rho - c_2L^2 \sqrt{\beta} \|w\|_{\psi_2} \sqrt{\log{N}} \frac{E_R}{\sqrt{N}}\right) \rho
\\
\geq & \left(c_1 \kappa_0^2 R - c_2L^2 \sqrt{\beta} \|w\|_{\psi_2} \sqrt{\log{N}} \frac{E_R}{\sqrt{N}}\right) R.
\end{align*} Therefore, if $N \gtrsim_{L,\kappa_0} E_R^2$ and
\begin{equation} \label{eq:E-R-in-proof}
E_R \leq c_3(L,\kappa_0) \frac{R}{\|w\|_{\psi_2}} \sqrt{\frac{N}{\beta\log N}},
\end{equation}
then $P_N {\cal L}_x > 0$ and $\hat{x} \not \in T_R$. Theorem~A follows from the definition of $r_2(\gamma)$ for a well chosen $\gamma$.


\section{Proof of Theorem B} \label{sec:proof-of-C}
Most of the work required for the proof of Theorem B has been carried out in Section \ref{sec:squared-loss}. A literally identical argument, in which one replaces the sets $T_{+,R}$ and $T_{-,R}$ with $T_{+,R}(x_0)$ and $T_{-,R}(x_0)$ may be used, leading to an analogous version of Theorem A, with the obvious modifications: the complexity parameter is $\max\{\ell(T_{+,R}(x_0)), \ell(T_{-,R}(x_0))\}$ for the right choice of $R$, and the probability estimate is $1-2\exp(-c_0 \min\{\ell^2(T_{+,R}(x_0)), \ell^2(T_{-,R}(x_0))\})-N^{-\beta+1}$.

All that remains to complete the proof of Theorem B is to analyze the structure of the local sets and identify the fixed points $r_0$ and $r_2$. A first step in that direction is the following:

\begin{Lemma} \label{localized-sets-large-norm}
There exist absolute constants $c_1$ and $c_2$ for which the following holds. For every $R>0$ and $\|x_0\|_2 \geq \sqrt{R}/4$,
\begin{description}
\item{1.} If $\|x_0\|_2 \min\{\|x-x_0\|_2, \|x+x_0\|_2 \} \geq R$ then $\|x-x_0\|_2 \|x+x_0\|_2 \geq c_1 R$.
\item{2.}  If $\|x-x_0\|_2 \|x+x_0\|_2 \geq R$ then $\|x_0\|_2 \min\{\|x-x_0\|_2, \|x+x_0\|_2 \} \geq c_2R$.
\end{description}
Moreover, if $\|x_0\|_2 \leq \sqrt{R}/4$ then $\|x-x_0\|_2\|x+x_0\|_2 \geq R$ if and only if $\|x\|_2 \gtrsim \sqrt{R}$.

\end{Lemma}

\proof
Assume without loss of generality that $\|x-x_0\|_2 \leq \|x+x_0\|_2$.

If $\|x-x_0\|_2 \leq \|x_0\|_2$ then
\begin{equation*}
\|x_0\|_2 \leq  2\|x_0\|_2 - \|x-x_0\|_2 \leq \|x+x_0\|_2 \leq \|x-x_0\|_2 + 2\|x_0\|_2 \leq 3\|x_0\|_2.
\end{equation*}
Hence, $\|x_0\|_2 \sim \|x+x_0\|_2$, and
$$
\|x_0\|_2 \min\{\|x-x_0\|_2,\|x+x_0\|_2\} \sim \|x-x_0\|_2 \|x+x_0\|_2.
$$

Otherwise, $\|x-x_0\|_2 > \|x_0\|_2$.

If, in addition,
$$
\|x_0\|_2 \geq (\|x-x_0\|_2 \|x+x_0\|_2)^{1/2}/4,
$$
then
$$
4 \|x_0\|_2 \geq (\|x-x_0\|_2 \|x+x_0\|_2)^{1/2} \geq \|x_0\|_2^{1/2} \|x+x_0\|_2^{1/2},
$$
and thus $\|x+x_0\|_2 \leq 16\|x_0\|_2$. Since $\|x_0\|_2 < \|x-x_0\|_2 \leq \|x+x_0\|_2$, it follows that $\|x+x_0\|_2 \sim \|x-x_0\|_2 \sim \|x_0\|_2$, and again,
$$
\|x_0\|_2 \min\{\|x-x_0\|_2,\|x+x_0\|_2\} \sim \|x-x_0\|_2 \|x+x_0\|_2.
$$
Therefore, the final case, and the only one in which there is no point-wise equivalence between $\|x-x_0\|_2 \|x+x_0\|_2$ and $\|x_0\|_2 \min\{\|x-x_0\|_2, \|x+x_0\|_2\}$,  is when $\min\{\|x-x_0\|_2,\|x+x_0\|_2\} \geq \|x_0\|_2$ and
$\|x_0\|_2 \leq (\|x-x_0\|_2 \|x+x_0\|_2)^{1/2}/4$. In that case, if $\|x_0\|_2 \geq \sqrt{R}/4$
then
$$
\|x_0\|_2 \min\{\|x-x_0\|_2, \|x+x_0\|_2\} \geq \|x_0\|_2^2 \geq R/16,
$$
and
$$
\|x-x_0\|_2 \|x+x_0\|_2 \geq 4\|x_0\|_2^2 \geq R/4,
$$
from which the first part of the claim follows immediately.

For the second one, observe that
$$
\|x\|_2^2 - 2\|x_0\|_2 \|x\|_2 \leq \|x-x_0\|_2 \|x+x_0\|_2 \leq \|x\|_2^2 +2 \|x\|_2 \|x_0\|_2 + \|x_0\|_2^2,
$$
and if $\|x_0\|_2 \leq \sqrt{R}/4$, the equivalence is evident.
\endproof

In view of Lemma \ref{localized-sets-large-norm}, the way the product $\|x-x_0\|_2\|x+x_0\|_2$  relates to $\min\{\|x-x_0\|_2, \|x+x_0\|_2\}$ depends on $\|x_0\|_2$. If $\|x_0\|_2 \geq \sqrt{R}/4$, then
$$
\{x \in T : \|x-x_0\|_2 \|x+x_0\|_2 \leq R\} \subset \{x \in T: \min\{\|x-x_0\|_2,\|x+x_0\|_2\} \leq c_1 R/\|x_0\|_2\},
$$
and if $\|x_0\|_2 \leq \sqrt{R}/4$,
$$
\{x \in T : \|x-x_0\|_2 \|x+x_0\|_2 \leq R\} \subset \{x \in T: \|x\|_2 \leq c_1 \sqrt{R}\},
$$
for a suitable absolute constant $c_1$.

When $T$ is convex and centrally-symmetric, the corresponding complexity parameter - the gaussian average of $T_{+,R}(x_0)=T_{-,R}(x_0)$ is
\begin{equation*}
E_R(x_0) \lesssim
\begin{cases}
\frac{\|x_0\|_2}{R} \cdot \ell(2T \cap (c_1R/\|x_0\|_2) B_2^n) & \mbox{if} \ \ \|x_0\|_2 \geq \sqrt{R},
\\
\\
\frac{1}{\sqrt{R}} \ell(2T \cap c_1\sqrt{R}B_2^n) & \mbox{if} \ \ \|x_0\| < \sqrt{R}.
\end{cases}
\end{equation*}
The fixed point conditions appearing in Theorem A now become
\begin{equation} \label{eq:r-0-thm-C}
r_0=\inf\{R: E_R(x_0) \leq c_2\sqrt{N}\}
\end{equation}
and
\begin{equation} \label{eq:r-2-thm-C}
r_2(\gamma) = \inf\{R: E_R(x_0) \leq \gamma \sqrt{N} R\},
\end{equation}
where one selects the slightly suboptimal $\gamma= c_2/ \sigma \sqrt{\log N}$. The assertion of Theorem A is that with high probability, ERM produces $\hat{x}$ for which
$$
\|\hat{x}-x_0\|_2 \|\hat{x}+x_0\|_2 \leq \max\{r_2(\gamma),r_0\}.
$$
If $\|x_0\|_2 \geq \sqrt{R}$, the fixed-point condition \eqref{eq:r-0-thm-C} is
\begin{equation} \label{eq:target-dependent-r-0}
\ell(2T \cap (c_1R/\|x_0\|_2) B_2^n) \leq c_3 \left(\frac{R}{\|x_0\|_2}\right) \sqrt{N}
\end{equation}
while \eqref{eq:r-2-thm-C} is,
\begin{equation} \label{eq:normalized-fixed-point}
\frac{\|x_0\|_2}{R} \ell(2T \cap (c_1R/\|x_0\|_2) B_2^n) \leq (c_4 / \sigma \sqrt{\log N}) \cdot \sqrt{N} R.
\end{equation}
Recall that
$$
r_N^*(Q) = \inf\{r>0: \ell(T \cap r B_2^n) \leq Q r \sqrt{N}\},
$$
and
$$
s_N^*(\eta)=\inf\{s>0 : \ell(T \cap s B_2^n) \leq \eta s^2 \sqrt{N}\}.
$$
Therefore, it is straightforward to verify that
$$
r_0=2\|x_0\|_2 r_N^*(c_3) \ \ {\rm and} \ \ r_2\big(c_2/(\sigma\sqrt{\log N})\big)=2\|x_0\|_2 s_N^*(c_4\|x_0\|_2/\sigma \sqrt{\log{N}}).
$$
Setting $R = 2\|x_0\|_2 \max\{r_N^*(c_3),s_N^*(c_4\|x_0\|_2/\sigma \sqrt{\log{N}})\}$, it remains to ensure that $\|x_0\|_2^2 \geq R$; that is,
\begin{equation} \label{eq:cond-x0}
2\max\{s_N^*(c_4\|x_0\|_2/\sigma \sqrt{\log N}),r_N^*(c_3)\} \leq \|x_0\|_2.
\end{equation}
Observe that if
\begin{equation} \label{eq:s-N-vs-r-N}
r_N^*(c_3) \leq \frac{c_3\sigma}{c_4\|x_0\|_2} \sqrt{\log N},
\end{equation}
then $r_N^*(c_3) \leq s_N^*(c_4\|x_0\|_2/\sigma \sqrt{\log N})$. Indeed, applying the convexity of $T$, it is standard to verify that $r_N^*(Q)$ is attained and $r_N^*(Q) \leq \rho$ if and only if $\ell(T \cap \rho B_2^n) \leq Q \rho \sqrt{N}$ -- with a similar statement for $s_N^*$ (see, e.g., the discussion in \cite{Lec-Men}). Therefore, $s_N^*(\eta) \geq r_N^*(Q)$ if and only if $\ell(T \cap r_N^*(Q)) \geq \eta (r_N^*(Q))^2 \sqrt{N}$. The latter is evident because $\ell(T \cap r_N^*(Q))=Qr_N^*(Q)\sqrt{N}$ and recalling that $Q=c_3$ and $\eta=c_4\|x_0\|_2/\sigma \sqrt{\log N}$.

Under (\ref{eq:s-N-vs-r-N}), an  assumption which has been made in the formulation of Theorem B, \eqref{eq:cond-x0} becomes $2s_N^*(c_4\|x_0\|_2/\sigma \sqrt{\log N}) \leq \|x_0\|_2$ and,  by the definition of $s_N^*$, this is the case if and only if
$$
\ell(T \cap \|x_0\|_2 B_2^n) \leq \frac{c_4 \|x_0\|_2}{\sigma \sqrt{\log N}} \cdot \|x_0\|_2^2 \sqrt{N};
$$
that is simply when $\|x_0\|_2 \geq v_N^*(\zeta)$ for $\zeta=c_4/\sigma \sqrt{\log N}$.

Hence, by Theorem A, combined with Lemma \ref{localized-sets-large-norm}, it follows that with high probability,
$$
\min\{\|\hat{x}-x_0\|_2,\|\hat{x}+x_0\|_2\} \leq 2s_N^*(c_4\|x_0\|_2/\sigma \sqrt{\log N}).
$$
The other cases, when either $\|x_0\|_2$ is `small', or when $r_0$ dominates $r_2$ are treated is a similar fashion, and are omitted.
\endproof

\section{Minimax lower bounds}
\label{sec:minimax}
In this section we study the optimality of ERM as a phase retrieval procedure, in the minimax sense. The estimate obtained here is based on the maximal cardinality of separated subsets of the class with respect to the $L_2(\mu)$ norm.
\begin{Definition} \label{def:packing}
Let $B$ be the unit ball in a normed space. For any subset $A$ of the space, let $M(A,rB)$ be the maximal cardinality of a subset of $A$ that is $r$-separated with respect to the norm associated with $B$.
\end{Definition}
Observe that if $M(A,rB) \geq L$ there are $x_1,...,x_L \in A$ for which the sets $x_i + (r/3)B$ are disjoint. A similar statement is true in the reverse direction.
\vskip0.3cm
Let $F$ be a class of functions on $(\Omega,\mu)$ and let $a$ be distributed according to $\mu$. For $f_0 \in F$ and a centred gaussian variable $w$, which has variance $\sigma$ and is independent of $a$, consider the gaussian regression model
\begin{equation}
  \label{eq:Gauss-model}
  y=f_0(a)+w.
\end{equation}
Any procedure that performs well in the minimax sense, must do so for any choice of $f_0\in F$  in \eqref{eq:Gauss-model}.

Following \cite{Lec-Men}, there are two possible sources of `statistical complexity' that influence the error rate of gaussian regression in $F$.
\begin{description}
\item{1.} Firstly, that there are functions in $F$ that, despite being far away from $f_{0}$, still satisfy $f_{0}(a_i)=f(a_i)$ for every $1 \leq i \leq N$, and thus are indistinguishable from $f_{0}$ on the data.

    This statistical complexity is independent of the noise, and for every $f_0 \in F$ and $\bA=(a_i)_{i=1}^N$, it is captured by the $L_2(\mu)$ diameter of the set
    $$
    K(f_0,\bA) = \{f \in F : (f(a_i))_{i=1}^N =(f_0(a_i))_{i=1}^N\},
    $$
which is denoted by $d_N^*(\bA)$.
\item{2.} Secondly, that the set $(F-f_{0})\cap r D =\{f-f_{0}: f\in F, \|f-f_{0}\|_{L_2} \leq r\}$ is `rich enough' at a scale that is proportional to its $L_2(\mu)$ diameter $r$.

    The richness of the set is measured using the cardinality of a maximal $L_2(\mu)$-separated set. To that end, let $D$ be the unit ball in $L_2(\mu)$, set
\begin{equation*}
C(r,\theta_0)=\sup_{f_0\in F}r\log^{1/2}M(F\cap(f_0+ \theta_0 r D),rD)
\end{equation*}
and put
\begin{equation}
  \label{eq:fixed-point}
  q^*_N(\eta)=\inf\big\{r>0: C(r,\theta_0)\leq \eta r^2\sqrt{N}\big\}.
\end{equation}

\end{description}
\begin{Theorem} \label{thm:LM-minimax} \cite{Lec-Men}
For every $f_0 \in F$ let $\Pro_{f_0}^{\otimes N}$ be the probability measure that generates samples $(a_i,y_i)_{i=1}^N$ according to \eqref{eq:Gauss-model}. For every $\theta_0 \geq 2$ there exists a constant $\theta_1>0$ for which
\begin{equation}\label{eq:minimax-bound}
  \inf_{\hat f}\sup_{f_0\in F}\Pro_{f_0}^{\otimes N}\Big(\norm{f_0-\hat f}_{2}\geq \max\{q_N^*(\theta_1/\sigma), (d_N^*(\bA)/4)\}\Big)\geq 1/5
\end{equation}
where $\inf_{\hat f}$ is the infimum over all possible estimators constructed using the given data.
\end{Theorem}
Earlier versions of this minimax bound may be found in \cite{MR2724359}, \cite{yang_barron99} and \cite{Bir83}.
\vskip0.3cm

To apply this general principle to the phase recovery problem, note that the regression function is  $f_0(x)=\inr{x_0,x}^2:=f_{x_0}(x)$ for some unknown vector $x_0\in T\subset\R^n$, while the estimators are $\hat f=\inr{\hat x,\cdot}^2$. Also, observe that for every $x_1, x_2 \in T$,
\begin{equation*}
  \norm{f_{x_0}-f_{x_1}}_{L^2(\mu)}^2=\E\big(\inr{x_0,a}^2-\inr{x_1,a}^2\big)^2=\E\inr{x_0-x_1,a}^2\inr{x_0+x_1,a}^2
\end{equation*} and  therefore, one has to identify the $L_2$ structure of the set
$$
F-f_{x_0}=\{\inr{x,\cdot}^2 - \inr{x_0,\cdot}^2 : x \in T\}.
$$
To obtain the desired bound, it suffices to assume the following:
\begin{Assumption}
  \label{ass:ass-isomorphy}
  There exist  constants $C_1$ and $C_2$ for which,
for every $s,t \in \R^n$,
\begin{equation*}
  C_1^2 \|s-t\|_2^2\|s+t\|_2^2\leq \E\inr{s-t,a}^2\inr{s+t,a}^2\leq C_2^2 \|s-t\|_2^2\|s+t\|_2^2.
\end{equation*}
\end{Assumption}

It is straightforward to verify that if $a$ is an $L$-subgaussian vector on $\R^n$ that satisfies Assumption \ref{ass:small-ball}, then it automatically satisfies Assumption \ref{ass:ass-isomorphy}.

The norm $\norm{x_0}_2$ plays a central role in the analysis of the rates of convergence of the ERM in phase recovery. Therefore, the minimax lower bounds presented here are not only for the entire model $T$ but for every shell $V_0=T\cap R_0 S^{n-1}$. A minimax lower bound over $T$ follows by taking the supremum over all possible choices of $R_0$.

To apply Theorem \ref{thm:LM-minimax}, observe that by Assumption~\ref{ass:ass-isomorphy}, for every $u,v \in T$,
$$
C_1 \|u-v\|_2 \|u+v\|_2 \leq \|f_{v}-f_u\|_{L_2} \leq C_2 \|u-v\|_2 \|u+v\|_2.
$$
Fix $R_0>0$ and consider $V_0=T\cap R_0S^{n-1}$. Clearly, for every $r>0$ and every $x_0 \in V_0$,
\begin{equation}
\label{eq:equiv-sets}
\Big\{u \in V_0 : \|u-x_0\|_2 \|u+x_0\|_2 \leq \frac{\theta_0 r}{C_2}\Big\} \subset \{u \in V_0: f_u \in F^\prime \cap (f_{x_0}+\theta_0 r D) \}
\end{equation}where $F^\prime=\{f_u:u\in V_0\}$.

Fix $\theta_0>2$ to be named later, and let $\theta_1$ be as in Theorem \ref{thm:LM-minimax}. If there are $x_0 \in V_0$ and $\{x_1,...,x_M\} \subset V_0$ that satisfy
\begin{description}
\item{1.} $\|x_i-x_0\|_2 \|x_i+x_0\|_2 \leq \theta_0 r/C_2$,
\item{2.} for every $1 \leq i < j \leq M$, $\|x_i-x_j\|_2 \|x_i+x_j\|_2 \geq r/C_1$, and
\item{3.} $ \log M > N (\theta_1 r/\sigma)^2$,
\end{description}
then $\sup_{f_0\in F^\prime}r\log^{1/2}M(F^\prime\cap(f_0+\theta_0 r D),rD) > \theta_1 r \sqrt{N}/\sigma$, and the best possible rate in phase recovery in $V_0$ is larger than $r$.

Fix $x_0 \in V_0$ and $r>0$, and let $R=r/C_2$. We will present two different estimates, based on $R_0=\norm{x_0}_2$, the `location' of $x_0$.
\vskip0.5cm
\noindent{\bf Centre of `small norm'.} Recall that $\theta_0>2$ and assume first that $R_0=\|x_0\|_2 \leq  \sqrt{\theta_0 R}/4$. Note that
$$
V_0 \cap (\sqrt{R}/8)B_2^n \subset \Big\{u \in V_0 : \|u-x_0\|_2 \|u+x_0\|_2 \leq \frac{\theta_0 r}{C_2}\Big\},
$$
and thus it suffices to constructed a separated set in $V_0 \cap (\sqrt{R}/8)B_2^n$.

Set $x_1,...,x_L$ to be a maximal $c_3\sqrt{R}$-separated subset of $V_0 \cap (\sqrt{R}/8)B_2^n$ for a constant $c_3$ that depends only on $C_1$ and $C_2$ and which will be specified later; thus, $L=M(V_0 \cap (\sqrt{R}/8)B_2^n,c_3\sqrt{R} B_2^n)$.
\begin{Lemma} \label{lemma:small-norm-set}
There is a subset $I \subset \{1,...,L\}$ of cardinality $M \geq L/2 -1$ for which $(x_i)_{i \in I}$ satisfies 1. and 2..
\end{Lemma}
\proof
Since $x_i \in (\sqrt{R}/8)B_2^n$ and $\|x_0\|_2 \leq \sqrt{\theta_0 R}/4$,
$$
\|x_i - x_0\|_2 \|x_i + x_0\|_2 \leq  ((\sqrt{\theta_0}+1)\sqrt{R}/4)^2 \leq \theta_0 R=\theta_0 r/C_2,
$$
and thus 1. is satisfied for every $1 \leq i \leq L$.

To show that 2. holds for a large subset, it suffices to find $I \subset \{1,...,L\}$ of cardinality at least $L/2-1$ such that for every $i,j \in I$,
$$
\|x_i-x_j\|_2 \geq c_3 \sqrt{R}/2 \ \ {\rm and} \ \ \|x_i+x_j\|_2 \geq c_3\sqrt{R}/2
$$
for a well-chosen $c_3$.

To construct the subset, observe that if there are distinct integers $1 \leq i,j, k \leq L$ for which $\|x_i + x_j\|_2 < c_3\sqrt{R}/2$ and   $\|x_i + x_k\|_2 < c_3\sqrt{R}/2$, then $\|x_j - x_k\|_2 < c_3\sqrt{R}$, which is impossible. Therefore, for every $x_i$ there is at most a single index $j \in \{1,...,L\} \backslash \{i\}$ satisfying that $\|x_i+x_j\|_2 < c_3\sqrt{R}/2$. With this observation the set $I$ is constructed inductively.

Without loss of generality, assume that $I=\{1,...,M\}$ for $M \geq L/2 -1$. If $i \not = j$ and $1 \leq i,j \leq M$,
$$
\|x_i-x_j\|_2 \|x_i+x_j\|_2 \geq c_3^2R/4 \geq 2r/C_1,
$$
for the right choice of $c_3$, and thus $(x_i)_{i=1}^M$ satisfies 2..
\endproof

\vskip0.5cm
\noindent{\bf Centre of `large norm'.} Next, assume that $R_0=\|x_0\|_2 \geq  \sqrt{\theta_0 R}/4$. By Lemma~\ref{localized-sets-large-norm}, there is an absolute constant $c_4 < 1/32$, for which, if $\|x_0\|_2 \geq \sqrt{\rho}/4$ and $\|x_0\|_2 \min\{\|x-x_0\|_2,\|x+x_0\|_2\} \leq c_4 \rho$, then $\|x-x_0\|_2 \|x+x_0\|_2 \leq \rho$.  Therefore, applied to the choice $\rho=\theta_0 R$,
\begin{align*}
& \left(V_0 \cap (x_0+(c_4\theta_0 R/\|x_0\|_2) B_2^n) \right) \cup \left(V_0 \cap (-x_0+(c_4\theta_0 R/\|x_0\|_2) B_2^n) \right)
\\
\subset & \left\{u \in V_0 : \|x_0-u\|_2 \|x_0+u\|_2 \leq \theta_0 R\right\},
\end{align*}
and it suffices to a find a separated set in the former.

Note that if $x \in V_0 \cap (x_0+(c_4\theta_0 R/\|x_0\|_2) B_2^n)$  then
$$
\|x\|_2 \geq \|x_0\|_2 - c_4\theta_0 R/\|x_0\|_2 \geq \|x_0\|_2/2 \geq \sqrt{\theta_0 R/4}/4
$$because one can choose $c_4\leq 1/32$, and $\|x\|_2 \leq 3\|x_0\|_2/2.$
Moreover, if $x_1,x_2 \in V_0 \cap (x_0+(c_4\theta_0 R/\|x_0\|_2) B_2^n)$, then
$$
\|x_1+x_2\|_2 \geq 2\|x_0\|_2 - 2c_4\theta_0 R/\|x_0\|_2 \geq \|x_0\|_2.
$$

Applying Lemma~\ref{localized-sets-large-norm}, there is an absolute constant $c_5$, for which, if
\begin{equation} \label{eq:separation-1}
\|x_0\|_2 \min \left\{\|x_i-x_j\|_2,\|x_i+x_j\|_2\right\} \geq \theta_0 R /4,
\end{equation}
then
$$
\|x_i-x_j\|_2 \|x_i+x_j\|_2 \geq c_5 \theta_0 R /4.
$$
Hence, if $x_1,...,x_M \in V_0 \cap (x_0+(c_4\theta_0 R/\|x_0\|_2) B_2^n)$ is $\theta_0 R/4\|x_0\|_2$-separated, then  \eqref{eq:separation-1} holds, and
$$
\|x_i-x_j\|_2 \|x_i+x_j\|_2 \geq c_5 \theta_0 R /4 \geq r/C_1,
$$
provided that $\theta_0$ is a sufficiently large constant that depends only on $C_1$ and $C_2$.


\begin{Corollary}\label{cor:minimax-pahse-reco1}
There exist absolute constants $\theta_1$, $c_1,c_2$ and $c_3$ that depend only on $C_1$ and $C_2$ and for which the following holds. Let $R_0>0$ and set $V_0=T\cap R_0S^{n-1}$. Define
\begin{equation*}
q^\prime(r)=\left\{
\begin{array}{cc}
\log M(V_0 \cap c_1 \sqrt{r}B_2^n,c_2\sqrt{r}B_2^n)   & \mbox{if} \ \ R_0 \leq c_3 \sqrt{r},
\\
 & \\
\sup_{x_0 \in V_0}\log M\left(V_0 \cap \left(x_0+c_1\left(\frac{r}{R_0}\right)B_2^n\right), c_2 \frac{r}{R_0}B_2^n\right)  & \mbox{if} \ \ R_0> c_3 \sqrt{r}
\end{array}\right.
\end{equation*}
If $q^\prime(r) \geq \theta_1 (r/\sigma)^2 N$, then the minimax rate in $V_0$ is larger than $r$.
\end{Corollary}

Now, a minimax lower bound for the risk $\min\{\norm{\tilde x-x_0}_2,\norm{\tilde x+x_0}_2\}$ for all shells $V_0$ (and therefore, for $T$) may be derived using Lemma~\ref{localized-sets-large-norm}. To that end, let $c_0$ be a large enough absolute constant and
\begin{equation}\label{eq:sudakov-complexity}
  C(R_0,r)=\sup_{x_0\in T:\norm{x_0}_2=R_0}r \log^{1/2} M\big((T\cap R_0S^{n-1})\cap(x_0+c_0 rB_2^n),rB_2^n).
\end{equation}
\begin{Definition}
  Fix $R_0>0$. For every $\alpha, \beta>0$ set
  \begin{equation*}
    q_N^*(\alpha)=\inf\big\{r>0:C(R_0,r)\leq \alpha r^2 \sqrt{N}\big\}
  \end{equation*}
  and put
  \begin{equation*}
    t_N^*(\beta)=\inf\big\{r>0:C(R_0,r)\leq \beta r^3 \sqrt{N}\big\}.
  \end{equation*}
\end{Definition}
Note that for any $c>0$, if $t_N^*(c)\leq1$ then $t_N^*(c)\geq q_N^*(c)$ and $q_N^*\big(cR_0/\sigma\big)\geq t_N^*(c/\sigma)$ if and only if $q_N^*\big(cR_0/\sigma\big)\geq R_0$.

\vskip0.7cm
\noindent{\bf Theorem C.} There exists an absolute constant $c_1$ for which the following holds. Let $R_0>0$.
\begin{enumerate}
\item If $R_0\geq t_N^*\big(c_1/\sigma\big)$, then for any procedure $\tilde x$, there exists $x_0\in T$ with $\norm{x_0}_2=R_0$ and for which, with probability at least $1/5$,
  \begin{equation*}
    \norm{\tilde x-x_0}_2\norm{\tilde x+x_0}_2\geq \norm{x_0}_2 q_N^*\Big(\frac{c_1\norm{x_0}_2}{\sigma}\Big)
  \end{equation*}
  and
  \begin{equation*}
   \min\{ \norm{\tilde x-x_0}_2,\norm{\tilde x+x_0}_2\}\geq  q_N^*\Big(\frac{c_1\norm{x_0}_2}{\sigma}\Big).
  \end{equation*}
\item If $R_0\leq t_N^*\big(c_1/\sigma\big)$ then for any procedure $\tilde x$ there exists $x_0\in T$ with $\norm{x_0}_2=R_0$ for which, with probability at least $1/5$,
  \begin{equation*}
    \norm{\tilde x-x_0}_2\norm{\tilde x+x_0}_2\geq  \Big(t_N^*\Big(\frac{c_1}{\sigma}\Big)\Big)^2
  \end{equation*}
  and
  \begin{equation*}
   \norm{\tilde x}_2, \norm{\tilde x-x_0}_2,\norm{\tilde x+x_0}_2\geq  t_N^*\Big(\frac{c_1}{\sigma}\Big).
  \end{equation*}
\end{enumerate}
\vskip0.7cm

Theorem~C is a general minimax bound, and although it seems strange at first glance, the parameters appearing in it are very close to those used in Theorem~C.  Following the same path as in \cite{Lec-Men}, let us show that Theorem~C and Theorem~C are almost sharp, under some mild structural assumptions on $T$.

First, recall Sudakov's inequality (see, e.g., \cite{LT}):
\begin{Theorem} \label{thm:sudakov}
If $W\subset\R^n$ and $\eps>0$ then
\begin{equation*}
  c\eps\log^{1/2} M(W,\eps B_2^n)\leq \ell(W),
\end{equation*}
where $c$ is an absolute constant.
\end{Theorem}

Fix $R_0>0$ set $V_0=T\cap R_0 B_2^n$, and put
$$
s_N^*=s_N^*\big(c_1R_0/(\sigma\sqrt{\log N})\big) \ \ {\rm and} \ \ v_N^*=v_N^*\big(c_1R_0/(\sigma\sqrt{\log N})\big).
$$
Assume that there is some $x_0\in T$, with the following properties:
\begin{description}
\item{1.} $\norm{x_0}_2=R$.
\item{2.} The localized sets $V_0\cap(x_0+s_N^*B_2^n)$ and $V_0\cap(x_0+v_N^*B_2^n)$ satisfy that
$$
\ell(V_0\cap (x_0+ s_N^*S^{n-1}))\sim \ell(T\cap s_N^* B_2^n)
$$
and that
$$
\ell(V_0\cap (x_0+v_N^*S^{n-1}))\sim \ell(T\cap v_N^* B_2^n),
$$
which is a mild assumption on the complexity structure of $T$.
\item{3.} Sudakov's inequality is sharp at the scales $s_N^*$ and $v_N^*$, namely,
\begin{equation}\label{eq:sudakov-sharp}
  s_N^*\log^{1/2}M(V_0\cap(x_0+c_0s_N^* B_2^n),s_N^* B_2^n)\sim \ell(V_0\cap(x_0+s_N^* B_2^n))
\end{equation}
and a similar assertion holds for $v_N^*$.
\end{description}
In such a case, the rates of convergence obtained in Theorem~B are minimax (up to the an extra $\sqrt{\log N}$ factor) thanks to Theorem~C.

\vskip0.3cm

When Sudakov's inequality is sharp as in \eqref{eq:sudakov-sharp}, we believe that ERM should be a minimax procedure in phase recovery, despite the logarithmic gap between Theorem~B an Theorem~C. A similar conclusion for linear regression was obtained in \cite{Lec-Men}.

Sudakov's inequality is sharp in many cases - most notably, when $T=B_1^n$, but not always. It is not sharp even for standard sets like the unit ball in $\ell_p^n$ for $1+(\log n)^{-1}<p<2$.


\section{Examples} \label{sec:examples}
Here, we will present two simple applications of the upper and lower bounds on the performance of ERM in phase recovery. Naturally, there are many other examples that follow in a similar way and that can be derived using very similar arguments. The choice of the examples has been made to illustrate the question of the optimality of Theorem~A, B and C, as well as an indication of the similarities and differences between phase recovery and linear regression. Since the estimate used in these examples are rather well known, some of the details will not be presented in full.

\subsection{Sparse vectors}

The first example we consider represents classes with a local complexity that remains unchanged, regardless of the choice of $x_0$.

Let $T=W_d$ be the set of $d$-sparse vectors in $\R^n$ (for some $d\leq N/4$) -- that is, vectors with at most $d$ non-zero coordinates. Clearly, for every $R>0$  $T_{+,R},T_{-,R} \subset W_{2d} \cap S^{n-1}$.
Also, for any $x_0\in T$ and any $I \subset \{1,...,n\}$ of cardinality $d$ that is disjoint of ${\rm supp}(x_0)$,
$$
(1/\sqrt{2}) S^I \subset \left\{ \frac{(x-x_0)_{i \in I}}{\|x-x_0\|_2} : x \in W_d\right\}, \left\{ \frac{(x+x_0)_{i \in I}}{\|x+x_0\|_2} : x \in W_d\right\},
$$
where $S^I$ is the unit sphere supported on the coordinates $I$.

With this observation, a straightforward argument shows that,
$$
\ell(T_{+,R}),\ell(T_{-,R}) \sim \sqrt{d \log (en/d)}.
$$
Applying Theorem~A, if follows that for $N\gtrsim d \log\big(en/d\big)$, with probability at least $1-2\exp(-c(d\log(en/d)))-N^{-\beta+1}$, ERM produces $\hat{x}$ that satisfies
\begin{equation} \label{eq:phase-sparse}
\|\hat{x}-x_0\|_2 \|\hat{x}+x_0\|_2 \lesssim_{\kappa_0,L,\beta} \sigma  \sqrt{\frac{d \log (en/d)}{N}}\sqrt{\log N}=(*).
\end{equation}
Moreover, with the same probability estimate, if $\norm{x_0}_2^2\gtrsim (*)$ then by Lemma \ref{localized-sets-large-norm},
\begin{equation}\label{eq:eq-large-x0-sparse}
\min\{\|\hat{x}-x_0\|_2, \|\hat{x}+x_0\|_2\} \lesssim_{\kappa_0,L,\beta} \frac{\sigma}{\norm{x_0}_2} \sqrt{\frac{d \log (en/d)}{N}}\sqrt{\log N}
\end{equation}
and if $\norm{x_0}_2^2\lesssim (*)$ then
\begin{equation}\label{eq:eq-small-x0-sparse}
\norm{\hat x}_2^2,\|\hat{x}-x_0\|_2^2, \|\hat{x}+x_0\|_2^2 \lesssim_{\kappa_0,L} \sigma \sqrt{\frac{d \log (en/d)}{N}}\sqrt{\log N}.
\end{equation}

In particular, when $\norm{x_0}_2$ is of the order of a constant, the rate of convergence in \eqref{eq:phase-sparse} and (\ref{eq:eq-large-x0-sparse}) is identical to the one obtained in \cite{Lec-Men} in the linear regression (up to a $\sqrt{\log N}$ term). In the latter, ERM achieves the minimax rate (with the same probability estimate) of
$$
\|\hat{x}-x_0\|_2 \lesssim_L \sigma \sqrt{\frac{d \log (en/d)}{N}}.
$$
Otherwise, when $\norm{x_0}_2$ is large, the rate of convergence of the ERM  in \eqref{eq:eq-large-x0-sparse}  is actually better than in linear regression, but, when $\norm{x_0}_2$ is small, it  is worse - deteriorating to the square root of the rate in linear regression (up to logarithmic terms).

When the noise level $\sigma$ tends to zero, the rates of convergence in linear regression and phase recovery tend to zero as well. In particular, exact reconstruction happens -- that is $\hat x=x_0$ in linear regression and $\hat x=x_0$ or $\hat x=-x_0$ in phase recovery -- when $N\gtrsim d\log\big(en/d\big)$.

For the lower bound, it is well known that $\log^{1/2} M(W_d\cap c_0 rB_2^n,rB_2^n)\sim \sqrt{d\log(en/d)}$ for every $r>0$ (and $c_0\geq2$). Combined with the results of the previous section, this suffices to show that the rate obtained in  Theorem~A is the minimax one (up to a $\sqrt{\log N}$ term in the ``large noise'' regime) 
and that the ERM is a minimax procedure for the phase retrieval problem when the signal $x_0$ is known to be $d$-sparse and $N\gtrsim d\log\big(en/d\big)$.

\subsection{The unit ball of $\ell_1^n$}
Consider the set $T=B_1^n$, the unit ball of $\ell_1^n$. Being convex and centrally symmetric, it is a natural example of a set with changing `local complexity' -- which becomes very large when $x_0$ is close to $0$. Moreover, it is an example in which one may obtain sharp estimates on $\ell(B_1^n \cap rB_2^n)$ at every scale. Indeed, one may show (see, for example, \cite{GLMP}) that
\begin{equation*}
  \ell\big(B_1^n\cap r B_2^n\big)\sim\left\{
    \begin{array}{cc}
      \sqrt{\log(enr^2)} & \mbox{ if } r^2n\geq1\\
      & \\
      r\sqrt{n}  & \mbox{ otherwise.}
    \end{array}
\right.
\end{equation*}
It follows that for $B_1^n$, one has
\begin{equation*}
  r_N^*(Q) \ \ \left\{
    \begin{array}{cc}
     \sim \Big(\frac{1}{Q^2 N}\log\Big(\frac{n}{Q^2 N}\Big)\Big)^{1/2} & \mbox{ if } n\geq C_0 Q^2 N\\
     & \\
\lesssim \frac{1}{N} & \mbox{ if } C_1Q^2N\leq n \leq C_0 Q^2 N\\
&\\
=0 & \mbox{ if } n\leq C_1 Q^2 N.
    \end{array}
\right.
\end{equation*}
where $C_0$ and $C_1$ are absolute constants. The only range in which this estimate is not sharp is when $n\sim Q^2N$, because in that range,  $r_N^*(Q)$ decays to zero very quickly. A more accurate estimate on $\ell(B_1^n\cap r B_2^n)$ can be performed when $n\sim Q^2N$ (see \cite{LM7}), but since it is not our main interest, we will not pursue it further, and only consider the cases $n\leq C_1 Q^2 N$ and $n\geq C_0 Q^2 N$.

A straightforward computation shows that the two other fixed points satisfy:
\begin{equation*}
  s_N^*(\eta)\sim\left\{
    \begin{array}{cc}
      \Big(\frac{1}{\eta^2 N}\log\Big(\frac{n^2}{\eta^2 N}\Big)\Big)^{1/4} & \mbox{ if } n\geq  \eta \sqrt{N}
      \\
      \\
  \sqrt{\frac{n}{\eta^2 N}} & \mbox{ if } n\leq  \eta \sqrt{N}
    \end{array}
\right.
\end{equation*}
and
\begin{equation*}
v_N^*(\zeta)\sim\left\{
    \begin{array}{cc}
      \Big(\frac{1}{\zeta^2 N}\log\Big(\frac{n^3}{\zeta^2 N}\Big)\Big)^{1/6} & \mbox{ if } n\geq  \zeta^{2/3} N^{1/3}
      \\
      \\
  \Big(\frac{n}{\zeta^2 N}\Big)^{1/4} & \mbox{ if } n\leq  \zeta^{2/3} N^{1/3}.
    \end{array}
\right.
\end{equation*}
The estimates above will be used to derive rates of convergence for the ERM $\hat x$ (for the squared loss) of the form
\begin{equation*}
  \min\{\norm{\hat x-x_0}_2,\norm{\hat x+x_0}_2\}\leq  rate.
\end{equation*}
Upper bounds on the rate of convergence ${\it rate}$ follow from Theorem~B, and hold with high probability as stated in there. For the sake of brevity, we will not present the probability estimates, but those can be easily derived from Theorem~B.

Thanks to Theorem~B, obtaining upper bounds on ${\it rate}$ involves the study of several different regimes, depending on $\norm{x_0}_2$, the noise level $\sigma$ and the way the number of observations $N$ compares with the dimension $n$.
\vskip0.3cm
\noindent{\bf The noise-free case: $\sigma=0$.} In this case, \textit{rate} is upper bounded by $r_N^*(Q)$, for $Q$ that is an absolute constant. In particular, when $n\geq C_0 Q^2 N$, the rate is less than
$$
\big(N^{-1}\log\big(n/N\big)\big)^{1/2}.
$$
At this point, it is natural to wonder whether there is a procedure that outperforms ERM in the noise-free case. The minimax lower bound $d_N^*(\bA)$ in Theorem~\ref{thm:LM-minimax} may be used to address this question, as no algorithm can do better than $d_N^*(\bA)/4$, with probability greater than $1/5$.

In the phase recovery problem  and using the notation of section~\ref{sec:minimax}, one has
\begin{align*}
  d_N^*(\bA)&=\sup \left\{\norm{f_x-f_{x_0}}_2: x\in B_1^n, \ f_x(a_i)=f_{x_0}(a_i) \ ,i=1,\ldots,N \right\}\\
&\sim \sup\left\{\norm{x-x_0}_2\norm{x+x_0}_2:\ x\in B_1^n, \ |\inr{a_i,x}|=|\inr{a_i,x_0}|, \ i=1,\ldots,N\right\}\\
&\gtrsim \inf_{L:\R^n\rightarrow\R^N}\sup\left\{\norm{x-x_0}_2\norm{x+x_0}_2: \ x\in B_1^n, \ L(x)=L(x_0)\right\}
\end{align*}
with an infimum taken over all linear operators $L:\R^n \to \R^N$.

By Lemma~\ref{localized-sets-large-norm}, for $x_0=(1/2,0,\ldots,0)\in B_1^n$ (in fact, any vector $x_0$ in $B_1^n$ for which $\norm{x_0}_2$ is a positive constant smaller than $1/2$ would do)
\begin{align*}
  d_N^*(\bA)&\gtrsim  \inf_{L:\R^n\rightarrow\R^N}\sup_{x\in B_1^n\cap({\rm ker}L-x_0)} \min\{\norm{x-x_0}_2,\norm{x+x_0}_2\}\\
&\gtrsim  \inf_{L:\R^n\rightarrow\R^N}\sup_{x,y\in B_1^n\cap{\rm ker}L} \norm{x-y}_2=c_N(B_1^n)
\end{align*}
which is the Gelfand $N$-width of $B_1^n$. By a result due to Garnaev and Gluskin (see \cite{MR759962}),
\begin{equation*}
  c_N(B_1^n)\sim \left\{
    \begin{array}{cc}
      \min\Big\{1,\sqrt{\frac{1}{N}\log\Big(\frac{en}{N}\Big)}\Big\} & \mbox{ if } N\leq n
      \\
      \\
0 & \mbox{ otherwise}.
    \end{array}
\right.
\end{equation*}
which is of the same order as $r_N^*$ (except when $n\sim N$, which is not treated here). Therefore, no algorithm can outperform ERM and ERM is a minimax procedure in this case.

Note that when $n\leq c_1 Q^2 N$, exact reconstruction of $x_0$ or $-x_0$ is possible and it can happens only in that case (i.e. $\sigma=0$ and $n\leq c_1 Q^2 N$) because of the minimax lower bound provided by $d_N^*(\bA)$.
\vskip0.3cm
\noindent{\bf The noisy case: $\sigma>0$.} According to Theorem~B, the rate of convergence \textit{rate} depends on $r_N^*=r_N^*(Q)$ for some absolute constant $Q$, $s_N^*=s_N^*(\eta)$ for $\eta=c_1 \norm{x_0}_2/(\sigma \sqrt{\log N})$ and on $v_N^*=v_N^*(\zeta)$ for $\zeta=c_1/(\sigma \sqrt{\log N})$. The outcome of Theorem~B is presented in Figure~1.
\begin{figure}[!h]\label{fig:rates-convergence}
\begin{center}
\begin{tabular}{|c|c|c|}
\hline
\textit{rate $\lesssim$} &   $\sigma/\norm{x_0}_2\leq c_0 r_N^*/\sqrt{\log N}$ & $\sigma/\norm{x_0}_2\geq c_0 r_N^*/\sqrt{\log N}$\\
\hline
$\norm{x_0}_2\leq v_N^*$ & $r_N^*$  &  $s_N^*$ \\
\hline
$\norm{x_0}_2\geq v_N^*$ & $r_N^*$   &  $v_N^*$\\
\hline
\end{tabular}
\end{center}
\caption{High probability bounds on the rate of convergence of the ERM $\hat x$ for the square loss in phase recovery: $\min\{\norm{\hat x-x_0}_2,\norm{\hat x+x_0}_2\}\leq rate$.}
\end{figure}

As the proof of all the estimates is similar, a detailed analysis is only presented when $\zeta^{2/3}N^{1/3}\leq \eta \sqrt{N}\leq C_1Q^2 N$, which is equivalent to
\begin{equation*}
  \Big(\frac{\sigma^2 \log N}{c_1 N}\Big)^{1/6}\leq \norm{x_0}_2\leq \frac{c_1Q^2 \sigma \sqrt{N \log N}}{c_1}.
\end{equation*}
The upper bounds on \textit{ rate } change according to the way the number of observations $N$ scales relative to $n$:
\begin{enumerate}
\item {\bf $n\geq C_0 Q^2 N$.} In this situation,  $r_N^*\sim\big(\log(n/N)/N\big)^{1/2}$. Therefore, if $\sigma/\norm{x_0}_2\lesssim \sqrt{\log(n/N)/(N\log N)}$  then $rate\leq \big(\log(n/N)/N\big)^{1/2}$, and if $\sigma/\norm{x_0}_2\gtrsim \sqrt{\log(n/N)/(N\log N)}$,
  \begin{equation}\label{eq:rate-conv-case1}
    rate\leq \left\{
      \begin{array}{cc}
        \Big(\frac{\sigma^2\log N}{\norm{x_0}_2^2 N}\log\Big(\frac{\sigma^2 n^2}{\norm{x_0}_2^2N}\Big)\Big)^{1/4} & \mbox{ if } \norm{x_0}_2\geq \Big(\frac{\sigma^2\log N}{ N}\log\Big(\frac{\sigma^2 n^3}{N}\Big)\Big)^{1/6}
        \\
        \\
\Big(\frac{\sigma^2\log N}{N}\log\Big(\frac{\sigma^2 n^3}{N}\Big)\Big)^{1/6} & \mbox{ otherwise}.
      \end{array}
\right.
  \end{equation}
\item {\bf $c_1\norm{x_0}_2/(\sigma \sqrt{\log N})\sqrt{N}\leq n\leq C_1Q^2 N$.} In that case $r_N^*=0$. In particular $\sigma/\norm{x_0}_2> c_0 r_N^*/\sqrt{\log N}$ and therefore, the rate is upper bounded as in (\ref{eq:rate-conv-case1}).
\item {\bf $\big(c_1/(\sigma \sqrt{\log N})\big)^{2/3}N^{1/3}\leq n\leq c_1\norm{x_0}_2/(\sigma \sqrt{\log N})\sqrt{N}$.} Again, in this case, $r_N^*=0$. Therefore, one is in the situation of the small \textit{signal-to-noise} ratio and
  \begin{equation*}
    rate\leq \left\{
      \begin{array}{cc}
    \frac{\sigma}{\norm{x_0}_2}\sqrt{\frac{n\log N}{N}}   & \mbox{ if } \norm{x_0}_2\geq \Big(\frac{\sigma^2\log N}{ N}\log\Big(\frac{\sigma^2 n^3}{N}\Big)\Big)^{1/6}
    \\
    \\
 \Big(\frac{\sigma^2\log N}{ N}\log\Big(\frac{\sigma^2 n^3}{N}\Big)\Big)^{1/6} & \mbox{ otherwise}.
      \end{array}
\right.
  \end{equation*}
\item {\bf $n\leq \big(c_1/(\sigma \sqrt{\log N})\big)^{2/3}N^{1/3}$.} Once again, $r_N^*=0$, and
  \begin{equation*}
    rate\leq \left\{
      \begin{array}{cc}
    \frac{\sigma}{\norm{x_0}_2}\sqrt{\frac{n\log N}{N}}   & \mbox{ if } \norm{x_0}_2\geq  \Big(\sigma \sqrt{\frac{n\log N}{N}}\Big)^{1/2}
    \\
    \\
 \Big(\sigma \sqrt{\frac{n\log N}{N}}\Big)^{1/2} & \mbox{ otherwise}.
      \end{array}
\right.
  \end{equation*}
\end{enumerate}

One may ask whether these estimates are optimal in the minimax sense, or perhaps there is another procedure that can outperform ERM. It appears that (up to an extra $\sqrt{\log N}$ factor), ERM is indeed optimal.

To see that, it is enough to apply Theorem~C and verify that Sudakov's inequality is sharp in the following sense: (see the discussion following Theorem~C): that if $\|x_0\|_1 \leq 1/2$, then for every $\eps<1/4$
\begin{equation*}
  \eps\log ^{1/2}M(B_1^n\cap(x_0+c_0\eps B_2^n),\eps B_2^n)\sim \ell\big(B_1^n\cap \eps B_2^n\big).
\end{equation*}
This fact is relatively straightforward to verify (see, e.g., Example 2 in \cite{Men-subgauss}).

Therefore, up to the extra $\sqrt{\log N}$ factor, which we believe is parasitic, ERM is a minimax phase-recovery procedure in $B_1^n$.

\begin{spacing}{0.9}
\begin{footnotesize}

\end{footnotesize}
\end{spacing}

\end{document}